\documentclass[a4paper,reqno]{amsart}

\usepackage{tgschola}

\usepackage{latexsym}
\usepackage[english]{babel}
\usepackage{fancyhdr}
\usepackage[mathscr]{eucal}
\usepackage{amsmath}
\usepackage{mathrsfs}
\usepackage{mathtools}
\usepackage{amsthm}
\usepackage{amsfonts}
\usepackage{amssymb}
\usepackage{amscd}
\usepackage{bbm}
\usepackage{graphicx}
\usepackage{graphics}
\usepackage{latexsym}
\usepackage{color}
\usepackage{pifont}
\usepackage{tikz}
\usepackage[normalem]{ulem}
\usepackage{changepage}

\usepackage{geometry}

\usepackage{makecell}

\usepackage{multirow}

\usepackage{tikz}
\usetikzlibrary{shapes, arrows.meta, decorations.markings, positioning, calc}

\theoremstyle{plain}

\theoremstyle{definition}

\newtheorem*{remark*}{Remark}

\numberwithin{equation}{section}

\begin{document}

\title[The Dynamical Landscape of Beggar-My-Neighbour]
{The Dynamical Landscape of Beggar-My-Neighbour: \\ Ultra-long Matches, Loops, and Infinite Matches}

\author[N.~Andorno]{Nicolas Andorno}
\address[N.~Andorno]{University of Trieste}
\email{nicolas.andorno@studenti.units.it}

\author[G.~Cernoia]{Giulio Cernoia}
\address[G.~Cernoia]{Department of Chemistry, Material, and Chemical Engineering \\  Polytechnic University of Milan, Milan \\ and European Space Agency, ESTEC \\ Utilisation and Enabled Science Team \\ Directorate of Human and Robotic Exploration \\ Noordwijk}
\email{giulio.cernoia@mail.polimi.it}

\author[S.~Duiz]{Simone Duiz}
\address[S.~Duiz]{University of Trieste}
\email{simone.duiz@studenti.units.it}

\author[A.~Michelangeli]{Alessandro Michelangeli}
\address[A.~Michelangeli]{Mathematics and Science Department, American University in Bulgaria, AUBG, 1 Georgi Izmirliev Sq., 2700 Blagoevgrad, \\
and Alexander von Humboldt Foundation, Bonn\\ and  Trieste Institute for Theoretical Quantum Technologies, TQT Trieste}
\email{amichelangeli@aubg.edu}

\date{\today}


%
%
%
%
%


\thanks{\emph{Acknowledgements.} For this work, A.M.~gratefully acknowledges support from the Italian National Institute for Higher Mathematics INdAM, the David Flanagan funds at the AUBG, American University Bulgaria, and the Alexander von Humboldt Foundation, Germany.}

\begin{abstract}
We present a rigorous mathematical and computational analysis of the deterministic card game \emph{Beggar-My-Neighbour}.
By establishing a formal state-space framework, we investigate the game's dynamical landscape, focussing on the dichotomy between terminating and non-terminating matches.
Extensive numerical simulations reveal that the distribution of finite match durations \emph{approximates} an exponential decay, with relevant deviations, confirming an emergent memory-less dynamics.
This statistical behaviour is further analysed in the context of ultra-long matches, where we identify characteristic multi-scale oscillatory patterns and entropic regimes.
Theoretically, we address the problem of backwards determinism, formalising the lack of injectivity of the trick function even within the set of reachable states.
Crucially, we contribute to the recent resolution of the long-standing question regarding the existence of infinite games.
We introduce an automated `Infinite Loop Factory' algorithm which, by implementing adaptive insertion strategies, proves effective in identifying non-terminating cycles with balanced initial deck configurations, thereby confirming the existence of non-terminating dynamics in standard and generalised settings of the game.
\end{abstract}

\maketitle

\tableofcontents

 \newpage

\section{Introduction}\label{sec:intro}

 \emph{Beggar-My-Neighbour} is a choice-free card game between two players, $A$ and $B$, where the cards played by each player follow deterministic rules based on the deck's initial setup; there are no actual free options for the players. The game is entirely predetermined as a sequence of actions from two sides.

 The game's rules, as detailed below, stipulate that players alternate revealing ordinary cards until one of the special cards is played, triggering a challenge phase in which the opponent must respond with a predetermined number of cards. Trick by trick, the winner collects the pile, and the match concludes when a player holds all the cards.

 Such a complete absence of player strategy -- as each match is choice-free and fully deterministic once the decks are prepared -- has not prevented Beggar-My-Neighbour from remaining a popular pastime since the mid-19th century, nor from attracting significant recent interest at both the amateur and academic levels.

 The drive behind the many recent challenges to find ultra-long matches, or even infinite loops, stems from the game's surprisingly rich dynamical landscape and the non-trivial complexity of the associated statistical and stochastic problems.

 These factors, in turn, motivated the present study, which analyses the game from multiple statistical viewpoints. A few precursor analyses in the literature and amateur contributions -- discussed in detail below -- have provided partial insights into specific aspects of the game, and this work aims to produce a unified outlook, supplemented by novel statistical analyses and findings.

 Here is the breakdown of our materials.

The game's rules are detailed in Section \ref{sec:game_rules}. We refer to the `standard' version of the game, which appears to enjoy the widest consensus, while briefly mentioning minority variants that are not considered in the following. We also consider a broad spectrum of scenarios depending on the number of ordinary and special cards in a deck; thus, our analysis focusses on a whole class of distinct settings, modelled on the original Beggar-My-Neighbour game.

A preliminary outlook on the game's probabilistic and statistical features is provided in Section \ref{sec:outlookanalysis}. This leads directly to Section \ref{sec:num}, which constitutes the first segment of our numerical analysis, based on a brute-force exploration of a portion of the vast space of initial configurations. The number of randomly selected distinct initial matches -- on the order of hundreds of millions -- provides a fairly clear picture of the distribution of match lengths. However, this remains a negligible fraction of the enormous number of possible initial decks, and is insufficient to detect the potentially extremely rare infinite matches.

Sections \ref{sec:abstr}-\ref{sec:fwd-bkw-determinism} constitute a theoretical block where the rigorous probabilistic formalism of the game is presented. Here, the focus is on: the total state space and the subspaces of initial, playable, and reachable configurations; the abstract action of the rule function that deterministically maps one trick to the next; the random variable of match duration, its distribution, and its expectation; and finally, the problem of the injectivity of the trick function -- a feature that, when present, guarantees the possibility of a uniquely determined backwards reconstruction of a game trajectory.

Section \ref{sec:et} reports the numerical analysis of ultra-long matches and discusses their characteristic patterns and distributions. Here we interpret the role of an approximate geometric distribution of the match duration, which explains the emergence of memory-less features in all considered indicators, such as the evolution of the card count, the average separation of special cards, and the position entropy of the decks.

The block of Sections \ref{sec:R1}-\ref{sec:inf-starting-loops} is entirely devoted to the problem of non-terminating matches. We begin with a thorough statistical analysis of the dynamics of two case studies: two recently identified infinite games, one found by us in a regime with a reduced number of special cards, and another in the original setting of the game. Phenomena of symmetry or imbalance in the evolution of the two players' decks are highlighted.

In Section \ref{sec:constructingInf}, we report on a very recent result from the literature regarding a suitably engineered construction of an infinite loop and of initial decks that evolve into such a loop after a transient phase. We reformulate that strategy in a manner we find more detailed and accessible, revealing that certain insightful ideas from that recent work, regarding the backward reconstruction of the loop, lay at the boundary between fully algorithmic analysis and successful heuristic intuition.

This motivates Section \ref{sec:loop_factory}, where we implement similar and new strategies for a more systematic algorithmic backward reconstruction. This successfully yields a variety of loops for multiple distinct settings of the game. Such loops are initially found to start from unbalanced decks, but a closer inspection in Section \ref{sec:inf-starting-loops} reveals that, in many of the loops we produced, balanced configurations are indeed reached, owing to the methodological innovations of our algorithmic search.

A final retrospective is provided in Section \ref{sec:retrosp-openpr}, covering both the main findings and perspectives of the present analysis, as well as a number of relevant open questions that naturally follow and to which we are committed.

These are challenges that may appear unique to the specific landscape of Beggar-My-Neighbour. Yet, their full understanding evidently requires advanced theoretical and numerical tools, and possibly the development of novel conceptual and practical approaches, at the intersection of mathematical logic, descriptive set theory, and combinatorial game theory.


%

\section{Game rules}\label{sec:game_rules}

\begin{figure}[!htbp]
    \centering

    \includegraphics[width=0.9\textwidth]{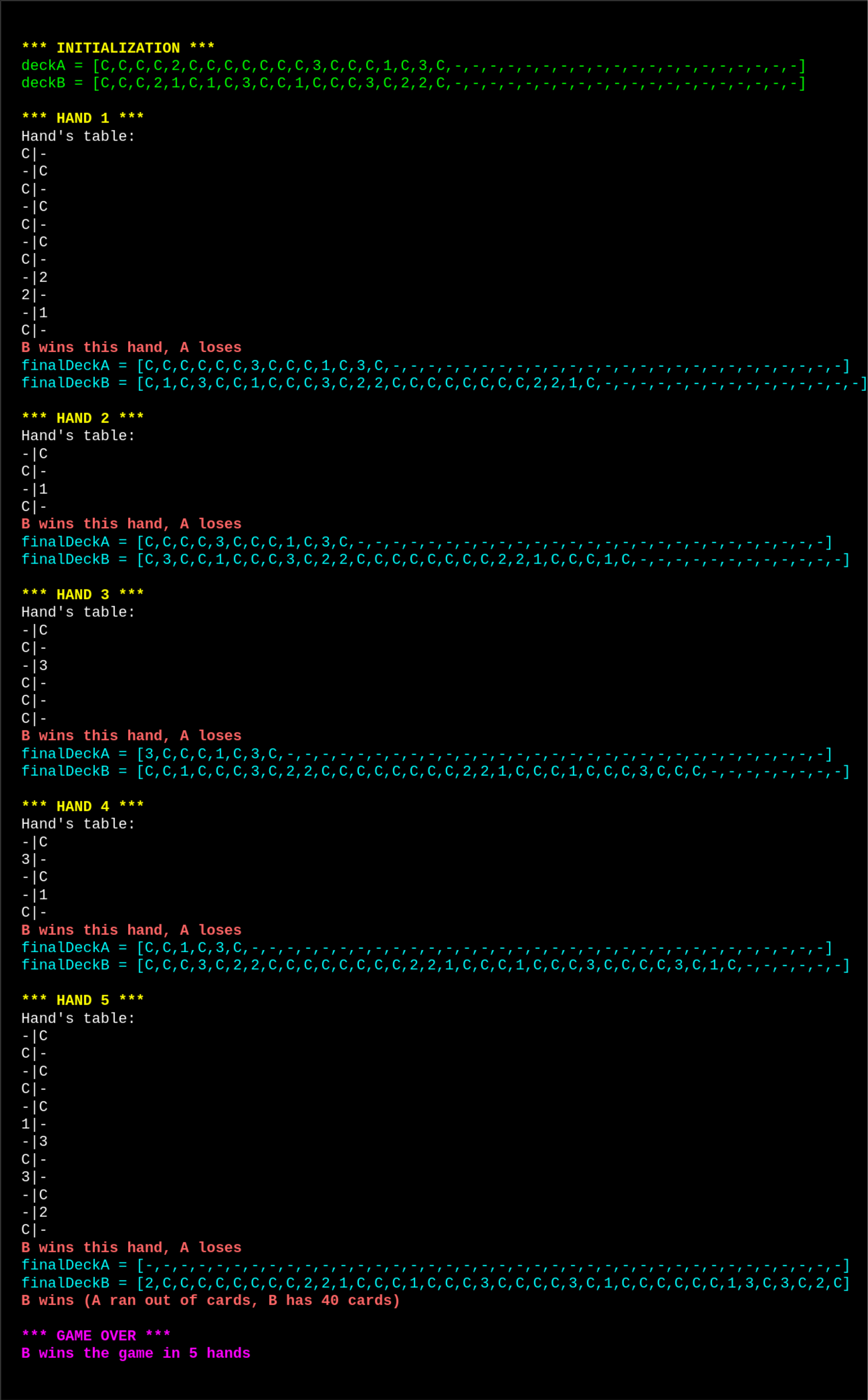}

    \caption{Example game with $N=40$ cards, and special cards up to maximum rank $\mathfrak{R}=3$, which is won in 5 tricks.}
    \label{fig:Exampleofgame}
\end{figure}


 The concrete version of Beggar-My-Neighbour considered here is a game played with a $N$-cards, French-type, four-suited deck, and $\mathfrak{R}$ special cards for each suit. Special cards are those ranked Ace through $\mathfrak{R}$ (e.g., if $\mathfrak{R}=3$, these are Aces, Twos, and Threes), and their number is $4\mathfrak{R}$. Ordinary cards (denoted by $C$) are all other cards, amounting to $N-4\mathfrak{R}$. We shall refer to this setting as the $(N,\mathfrak{R})$-setting of the game.

 At the beginning of each match, the deck is randomised by shuffling and then split evenly. Conventionally, the player starting the first trick of a match shall be labelled player $A$ for that match. For all subsequent tricks in the same match, the player who won the previous trick goes first. At the beginning, player $A$ receives the first $N/2$ cards face down, and player $B$ receives the other $N/2$ cards, also face down.

 A trick starts with the opening player flipping their top card face up. The other player then flips their top card face up on top of the first, and this back-and-forth continues as long as only ordinary cards are played.

If a special card, say of rank $r$, $r\in\{1,\dots,\mathfrak{R}\}$, is played, irrespective of its suit, the ongoing alternating card-flipping phase pauses. The player who revealed the special card is now designated the `challenger', and the other player becomes the `opponent'. The opponent must now attempt to respond, playing a maximum of $r$ consecutive cards from their deck. Should all of these $r$ cards be ordinary $C$-cards, the opponent's response is unsuccessful, and the trick concludes. The challenger wins this trick. If instead, at any point during their response, the opponent plays a special card of rank $r'$, the current challenge is immediately terminated. The match then continues with a new challenge, but with the roles swapped: the player who was just the challenger now becomes the opponent, and must respond with up to $r'$ cards from their deck, and so forth.

Once a trick is completed, the pile of cards on the table is collected, turned face down, and added to the bottom of the winner's remaining deck.

When a player runs out of cards (i.e., their deck is empty), they lose, and the match ends. The winner is the player who has accumulated all $N$ of the original cards. Figure \ref{fig:Exampleofgame} illustrates an example game in the $(40,3)$-setting.

 We remark that the traditional Beggar-My-Neighbour setting corresponds to the $(52,4)$-configuration, with special cards identified as Jack (J), Queen (Q), King (K), and Ace (A), which in our notation have ranks 1, 2, 3, and 4, respectively.

 The $(40,3)$-setting represents the Italian variant of the game, known, among others, by the names \emph{Cavacamisa}, \emph{Pataiola}, \emph{Pataja}, \emph{Pelagallina}, \emph{Rubamazzetto}, \emph{Scipacuore}, \emph{Stracciacamicia}, and \emph{Andogame}.

 For clarity, we also summarise the terminology that shall be used throughout:
 \begin{itemize}
  \item \emph{move}: the play of a single card from the top of a player's deck to the central pile;
  \item \emph{trick} or \emph{hand}: a complete sequence of moves that begins when the central pile is empty and ends the first time one player wins all cards from the pile and adds them to the bottom of their deck;
  \item \emph{match} (or just \emph{game}, if it is clear from the context):  a complete game that starts from the preparation of two initial half-decks, unfolds across a sequence of tricks where each trick begins with a move by the player who won the previous trick, and ends when a player wins a trick that leaves their opponent without any remaining cards.
 \end{itemize}

\section{Outlook analysis}\label{sec:outlookanalysis}

  Beggar-My-Neighbour evolves with dynamics solely determined by the occurrence of special cards.

  The number $\mathcal{N}_{\mathrm{init}}$ of special cards in a player's \emph{initial} $N/2$-cards deck follows a hypergeometric distribution $\mathcal{N}_{\mathrm{init}}\sim\textrm{Hypergeometric}(N,4\mathfrak{R},N/2)$ (total population of $N$ cards, successes of $4\mathfrak{R}$ cards, sample size of $N/2$ cards).

  Conditional to the value of $\mathcal{N}_{\mathrm{init}}$ the position of the first special card in each player's half-deck follows a negative hypergeometric distribution (sampling without replacement), which implies that the probability of finding the first special card deeper in the stack decreases faster than in the geometric case.

  The evolution of the players' decks in subsequent tricks, driven by challenge phases, distorts these initial distributions in a manner that renders closed-form analysis intractable.
Obtaining a probabilistic description of each trick, including the number of cards played and the trick's winner, is computationally and combinatorially explosive.

In turn, the number and distribution of special cards in each player's deck at each trick crucially affect the length of each match.

In one direction, a one-trick match is always possible, for example with players' decks consisting of all generic cards on top, and all special cards at the bottom.

On the other hand, the actual length of a match can be much longer. It is reasonable to argue that the larger the ratio $\mathfrak{R}/N$ the more likely the occurrence of challenges, with cards being consumed earlier, and with the net effect of making the match shorter. Conversely, a small $\mathfrak{R}/N$ ratio results in longer alternating card-flipping phases and shorter challenge phases, thereby producing matches with large amounts of tricks.

These expectations are sharply confirmed by the numerical analysis presented in Sections \ref{sec:num} and \ref{sec:et}. Yet, on top of such considerations, the length of a match is also non-trivially affected by the distribution of special cards.

Consequently, in light of the complexities discussed above, the general problem of the length of a match is a hard one.

This extends to the possibility of never-ending matches, namely matches with an infinite number of tricks.

The existence of such infinite trajectories -- cycles in the state space -- would fundamentally alter the statistical description of the game, implying that the expected game duration is formally infinite. Determining whether such cycles exist for a given $(N, \mathfrak{R})$ configuration has been a long-standing open problem in combinatorial game theory, until the recent identification of loops.

\section{Numerical analysis up to hundreds of millions of matches}\label{sec:num}

 The space of distinct initial deck configurations is vast: their number $\mathcal{C}(N,\mathfrak{R})$ amounts to
 \begin{equation}
  \mathcal{C}(N,\mathfrak{R})\:=\:\frac{N!}{\,(4!)^{\mathfrak{R}} (N-4\mathfrak{R})!\,}
 \end{equation}
 (with the convention of having 4 suits).
 For example, $\mathcal{C}(52,4)\sim 10^{21}$ and  $\mathcal{C}(40,3)\sim 10^{14}$, to mention two most typical settings. (It is worth remarking the somewhat counterintuitive pattern
 \[
  \mathcal{C}(52,2)\sim 10^{11}\,,\quad
    \mathcal{C}(52,3)\sim 10^{16}\,,\quad
      \mathcal{C}(52,4)\sim 10^{21}\,,\quad
        \mathcal{C}(52,5)\sim 10^{26}\,,\quad
          \mathcal{C}(52,6)\sim 10^{30}\,,\quad
 \]
 etc., that is, adding structure with additional types of special cards actually  increases rather than decreases the number of distinguishable arrangements.)

Brute force numerics, even when pushed to billions of distinct matches, barely sample the space of possible initial decks. Among past computational studies,
\begin{itemize}
  \item \cite{Paulhus_Beggar1999} declared simulation of $\sim 10^{9}$ distinct decks for the $N=52$, $\mathfrak{R}=4$ setting;
  \item \cite{Beggar-Collins-records} declared simulation of $\sim 10^{11}$ distinct decks for the $N=52$, $\mathfrak{R}=3$ setting;
  \item \cite{Beggar-Zanotto-codes} declared simulation of $\sim 10^{7}$ distinct decks for the $N=40$, $\mathfrak{R}=3$ setting.
\end{itemize}
Additionally, \cite{Casella2024} reported, based on personal correspondence, claimed simulations by others of the $(52,4)$-game for $10^{13}$ to $10^{15}$ matches.

 We present here evidence from simulating $10^7$-$10^8$ randomly selected distinct matches, in dedicated C scripts. Significantly, this amount of simulated matches proves not to be large enough to allow for the detection of the occurrence of infinite matches. It therefore provides a statistical analysis of finite matches.

\begin{figure}[!htbp]
    \begin{adjustwidth}{-2.4cm}{-3.7cm}
        \centering
        \begin{minipage}[b]{0.55\textwidth} 
            \centering
            \includegraphics[width=\linewidth]{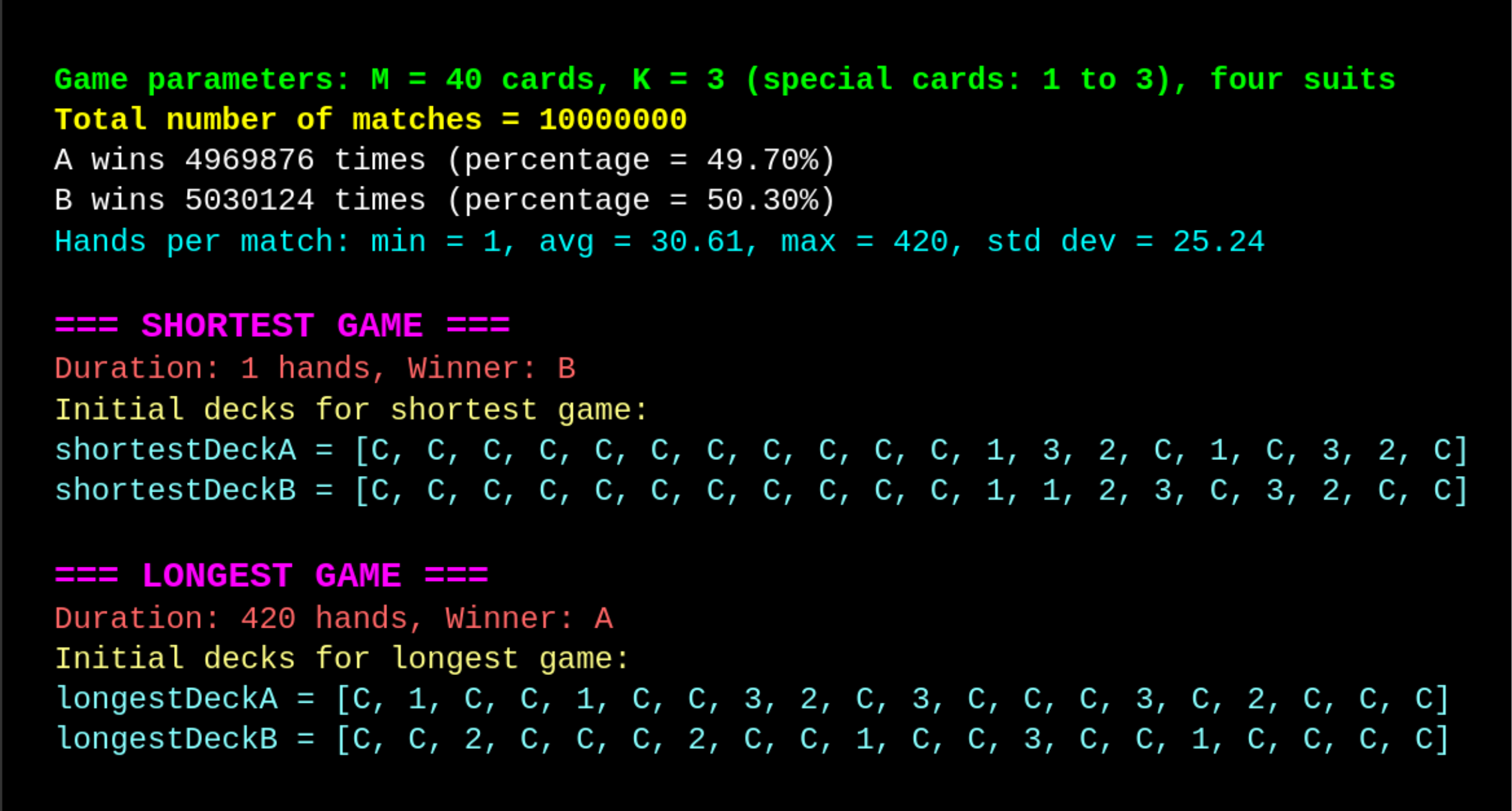} 
        \vspace{0.1cm}
        \end{minipage}
        \quad 
        \begin{minipage}[b]{0.75\textwidth} 
            \centering
            \includegraphics[width=\linewidth]{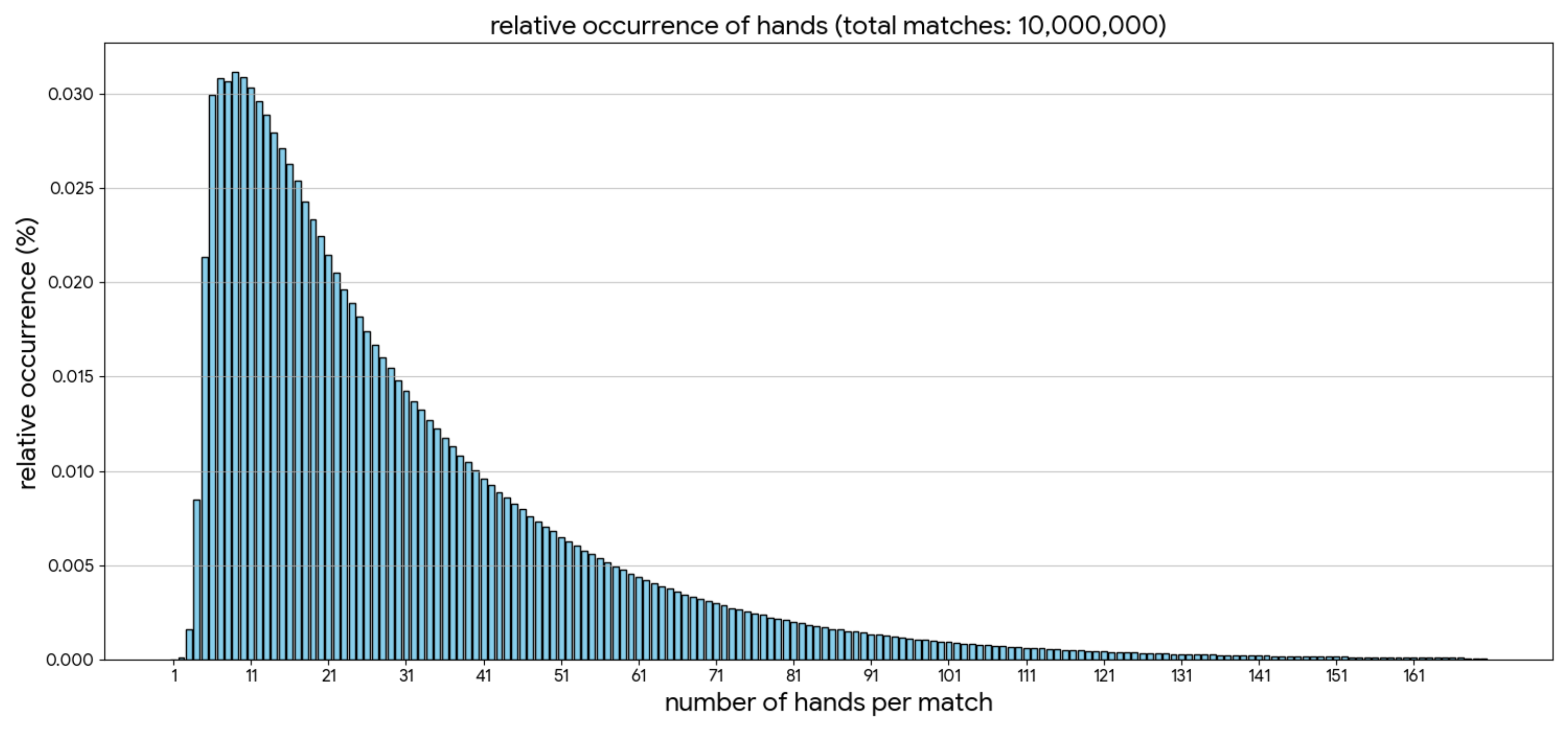} 
        \end{minipage}
        \caption{\small Left: Output of a numerical simulation of $10^7$ independent matches for the game with $N=40$ cards, and special cards' highest rank $\mathfrak{R}=3$. Right: frequency distribution of the matches versus their length in number of tricks  per match. \normalsize}
        \label{fig:stats_output_figure}
    \end{adjustwidth}
\end{figure}

Figure \ref{fig:stats_output_figure} displays a typical simulation output, demonstrating the possibility of one-trick matches as well as the occurrence of much longer matches. The frequency of played tricks per match exhibits a characteristic right-skewed profile.

\begin{figure}[!htbp]
    \centering

    \includegraphics[width=0.7\textwidth]{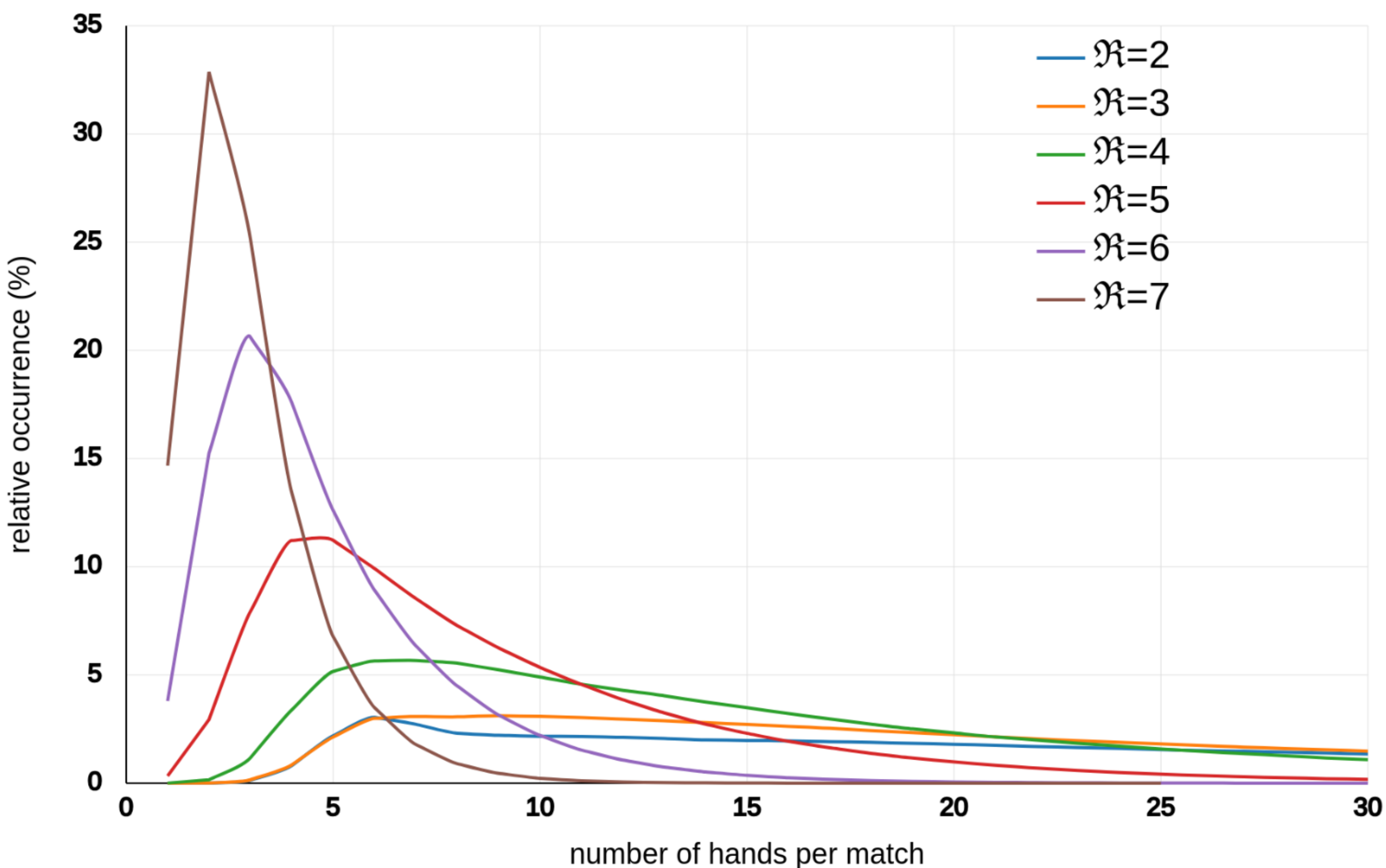}

    \caption{\small Frequency distribution of the number of tricks per match for each $(N,\mathfrak{R})$-setting game, with $N=40$ and $\mathfrak{R}\in\{2,3,4,5,6,7\}$. \normalsize}
    \label{fig:r_parameter_plot}
\end{figure}

 \begin{figure}[!htbp]
    \centering
    \includegraphics[width=0.8\textwidth]{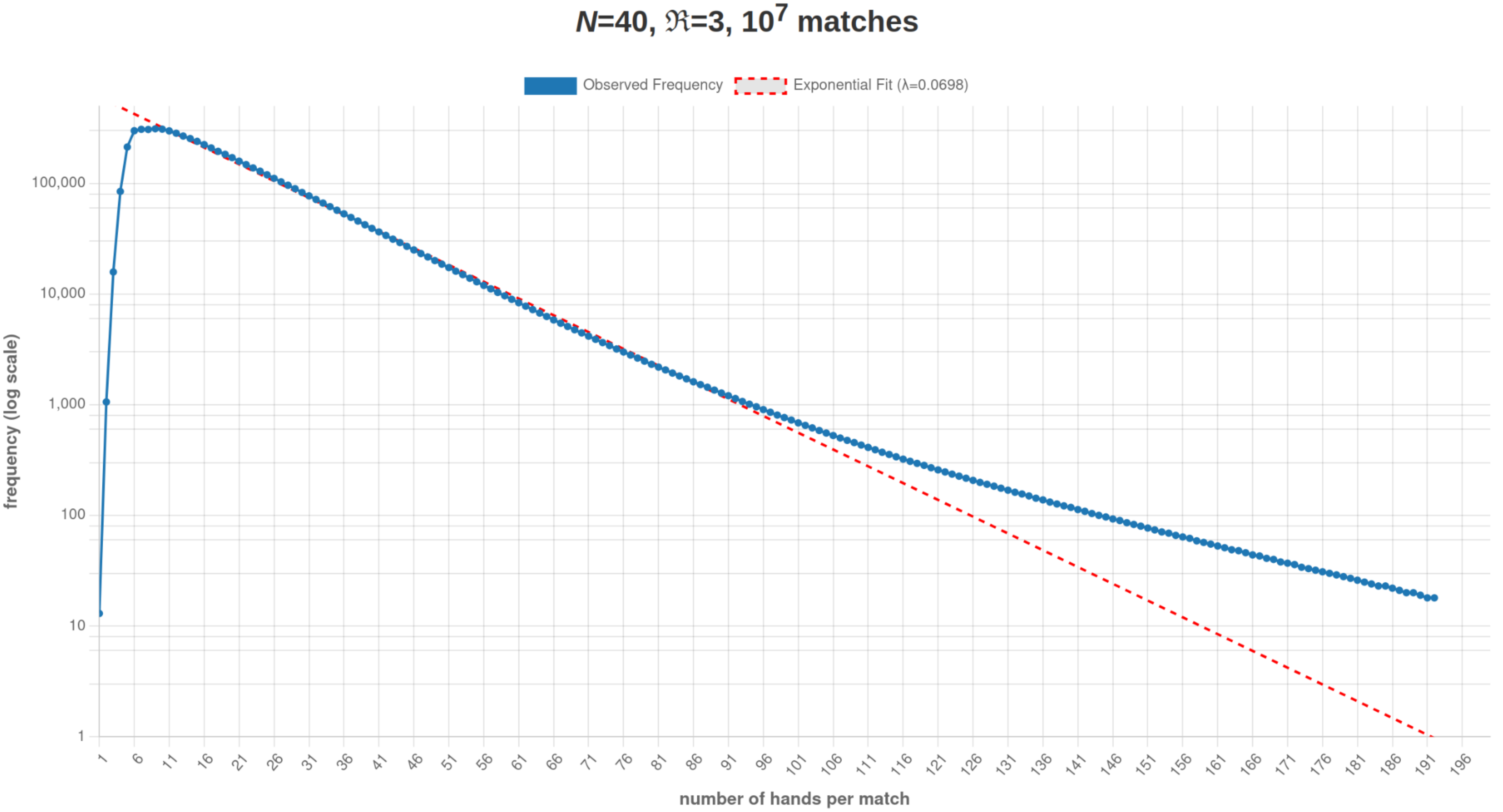}
    \caption{Semi-logarithmic plot of game duration frequency for the $(40,3)$-setting based on the matches simulation reported in Table \ref{tab:numanalrecap} and \ref{tab:numanalrecapN}. The approximately linear decay in the range $n \in [10, 100]$ is consistent with exponential distribution $\rho(n) \propto e^{-\lambda n}$ with $\lambda \approx 0.0394$. Strong deviation from linearity occurs for $n \leqslant 9$ (short games).}
    \label{fig:semilog_40_3}
\end{figure}

\begin{table}[ht!]
\centering
\begin{tabular}{||l||c|c|c|c|c|c|c||}
\hline
$N$ (total cards) & 40 & 40 & 40 & 40 & 40 & 40 & 40 \\
\hline
$\mathfrak{R}$ (highest special rank) & 2 & 2 & 3 & 4 & 5 & 6 & 7 \\
\hline
special cards & 1 to 2 & 1 to 2 & 1 to 3 & 1 to 4 & 1 to 5 & 1 to 6 & 1 to 7 \\
\hline
simulated matches & $10^8$ & $10^7$ & $10^7$ & $10^7$ & $10^7$ & $10^7$ & $10^7$ \\
\hline
min number of tricks per match & 1 & 1 & 1 & 1 & 1 & 1 & 1 \\
\hline
max number of tricks per match & 700 & 634 & 420 & 215 & 91 & 49 & 25 \\
\hline
most likely no.~of tricks per match & 6 & 6 & 9 & 7 & 5 & 3 & 2 \\
\hline
average number of tricks per match & 41.26 & 41.26 & 30.61 & 17.41 & 8.92 & 4.81 & 2.93 \\
\hline
standard deviation & 36.94 & 36.62 & 25.24 & 12.95 & 5.95 & 2.90 & 1.60 \\
\hline
variance-to-mean ratio & 33.07 & 33.11 & 21.12 & 9.85 & 4.01 & 1.78 & 0.89 \\
\hline
A wins (\%) & 49.78 & 49.76 & 49.70 & 49.66 & 49.61 & 49.39 & 48.65 \\
\hline
B wins (\%) & 50.22 & 50.24 & 50.30 & 50.34 & 50.39 & 50.61 & 51.35 \\
\hline
\end{tabular}
\vspace{0.2cm}
\caption{\small Summary of the numerical and statistical analysis in various $(40,\mathfrak{R})$-settings.\normalsize}
\label{tab:numanalrecap}
\end{table}

\begin{table}[ht!]
\centering
\begin{tabular}{||l||c|c|c|c|c|c|c||}
\hline
$N$ (total cards)   & 52  & 48  & 44 & 40 &36 & 32 & 28\\
\hline
$\mathfrak{R}$ (highest special rank)  & 3 & 3 & 3 & 3 & 3 & 3 & 3 \\
\hline
special cards  & 1 to 3 & 1 to 3 & 1 to 3 & 1 to 3 & 1 to 3 & 1 to 3 & 1 to 3 \\
\hline
simulated matches  & $10^7$ & $10^7$ & $10^7$ & $10^7$  & $10^7$ & $10^7$ & $10^7$\\
\hline
min number of tricks per match  & 2 & 1 & 1 & 1 & 1 & 1 & 1 \\
\hline
max number of tricks per match   & 706 & 645 & 513  & 420  & 326 & 284  & 178 \\
\hline
most likely no.~of tricks per match  & 15 & 12  & 10  & 9  & 6 & 6 & 5 \\
\hline
average number of tricks per match    & 51.16  & 44.20  & 37.32  & 30.61 & 24.14 & 18.07 & 12.56 \\
\hline
standard deviation    & 43.66  & 37.36  & 31.29  & 25.24 & 19.41 & 14.10 & 9.05 \\
\hline
variance-to-mean ratio   & 37.76 & 32.11  & 26.56  & 21.12 & 15.89 & 11.00 & 6.72 \\
\hline
A wins (\%)  & 49.81  & 49.77  & 49.71  & 49.70 & 49.62 & 49.53 & 49.43 \\
\hline
B wins (\%)  & 50.19  & 50.23  & 50.29  & 50.30  & 50.38 & 50.47 & 50.57 \\
\hline
\end{tabular}
\vspace{0.2cm}
\caption{\small Summary of the numerical and statistical analysis in various $(N,3)$-settings.\normalsize}
\label{tab:numanalrecapN}
\end{table}

The comprehensive simulation results from $10^7$-$10^8$ independent matches for each $(N,\mathfrak{R})$-setting game are presented in Table \ref{tab:numanalrecap} for $N=40$ and $\mathfrak{R}\in\{2,3,4,5,6,7\}$, and in Table \ref{tab:numanalrecapN} for $N\in\{28,32,36,40,44,48,52\}$ and $\mathfrak{R}=3$, in all cases with a four-suited deck. In both tables the ratio $4\mathfrak{R}/N$ of special vs total cards increases from left to right.

A clear pattern emerges, in which, as argued in Section \ref{sec:outlookanalysis}, as $\mathfrak{R}/N$ increases the average number of tricks per match decreases (and its standard deviation also decreases).

One-trick-only (minimum-length) matches were consistently detected (though never in the $(52,3)$-setting), and the maximum number of tricks per match substantially increases with decreasing $\mathfrak{R}$.

Interestingly, whereas players $A$ and $B$ achieve an approximate $50\%$-$50\%$ share of victories over a large number of matches, player $A$ (the one who makes the first move in the first trick of a match) consistently has a victory percentage just under $50\%$, which furthermore slightly decreases as $\mathfrak{R}$ increases. The opposite trend is observed for player $B$. Playing the first move in a match brings on average a tiny handicap.

The increase, as $\mathfrak{R}/N$ increases, of the right-skewed character of the distribution of the number of tricks per match is shown in Figure \ref{fig:r_parameter_plot}. The peaks in each $(N,\mathfrak{R})$-curve indicate the most likely number of tricks per match, as also reported in Tables \ref{tab:numanalrecap} and \ref{tab:numanalrecapN}.


 The pronounced right skewness of the distribution of match lengths is more accurately fitted by an exponential distribution, rather than a Poisson distribution. The latter should be excluded by standard variance-to-mean dispersion tests. In a genuine Poisson scenario, the mean and variance coincide. However, here the variance-to-mean ratio is significantly greater than 1 for special-vs-total-cards ratio $4\mathfrak{R}/N\leqslant 50\%$ (for concreteness: $2\leqslant\mathfrak{R}\leqslant 5$ when $N=40$), indicating over-dispersion and thus much more variability than Poisson predicts. This is likely due to clustering/contagion (events occurring in bursts) and non-independence of trick counts. Conversely, for $4\mathfrak{R}/N\geqslant 70\%$ (for concreteness: $\mathfrak{R}\geqslant 7$ when $N=40$), the ratio is less than 1, indicating under-dispersion (less variability than Poisson), which is likely caused by the regulation of extreme values and limited `opportunities' for trick counts.

 It is rather the exponential fit that turns out to be accurate, at least throughout a wide intermediate range of game durations, from a few units to hundreds, albeit with heavier tails (more long games). Semi-logarithmic plots of the reported simulations reveal this good agreement (Figure \ref{fig:semilog_40_3}). This is an observation already emerged in multiple amateur analyses (see, e.g., \cite{Beggar-Gentilini}), with which our data are well consistent. Recently, \cite{Casella2024} reproduced this observation also restricting to initial decks prepared with the same number of special cards: the empirical distribution of match lengths remains essentially unaltered.

 \section{Game analysis' abstract framework}\label{sec:abstr}

 In parallel to numerics and empirical study, a convenient abstract framework for the rigorous mathematical analysis of Beggar-My-Neighbour appears to have been elaborated first in \cite{Lakshtanov-Aleksenko-Beggar2013}. For the purposes of that work, the game's rules were substantially altered, as discussed here later on. In this Section the original formalism of \cite{Lakshtanov-Aleksenko-Beggar2013} is re-formulated and expanded so as to account for the correct version of the game.

Let $D$ denote the complete deck for the $(N,\mathfrak{R})$-setting of the game, and let $\widehat{D}$ be a re-shuffling of the deck, namely $D$ itself, equipped with the position ordering of its cards (thus, the various $\widehat{D}$'s are just permutations of each other).

Tracking the ordering is essential because
\begin{itemize}
\item the card position determines accessibility (players can only access their top card),
\item the sequence order determines all future play possibilities,
\item and state transitions depend critically on the exact ordering when cards are added to deck bottoms.
\end{itemize}

 The state space $\mathcal{S}\equiv\mathcal{S}_{N,\mathfrak{R}}$ of the game, meant as the collection of all possible ordered configurations of splittings of $D$ among two subjects $A$ and $B$, is then described as
\begin{equation}
\mathcal{S} = \bigcup_{\widehat{D} \in \text{Perm}(D)} \{(L,R) : L \text{ and } R \text{ are ordered sequences, } L \cdot R = \widehat{D}\}\,,
\end{equation}
where
\begin{itemize}
\item $L$ represents player $A$'s deck as an ordered sequence (top to bottom),
\item $R$ represents player $B$'s deck as an ordered sequence (top to bottom),
\item $L \cdot R$ denotes concatenation reconstructing the full deck $\widehat{D}$,
\item hence the sizes $|L|$ and $|R|$ of the two players' decks are integers in $\{0,1,\dots,N\}$ and add up to $|L|+|R|=N$.
\end{itemize}
The notation is to indicate that `$A$ plays sitting at the \emph{L}eft', and `$B$ at the \emph{R}ight'.

 In practice, only certain subspaces of $\mathcal{S}$ are relevant in the game (most importantly: starting-game decks ($|L|=|R|=\frac{N}{2}$), intermediate decks that are actually reachable along a match history, decks that can be actually produced by previous decks through one single trick). This will be discussed in the following.

Denote by $\mathcal{S}^*$, in the conventional Kleene star notation, the collection of all possible \emph{finite} sequences in $\mathcal{S}$. The match evolution over finitely many tricks is represented by a sequence $\sigma = (s_0, s_1, s_2, \ldots, s_{\texttt{end}})\in\mathcal{S}^*$ of states with $s_j = (L_j, R_j) \in \mathcal{S}$ denoting the decks configuration after completing the $j$-th trick. Here $s_0$ is an initial configuration at the beginning of the match, and $s_{\texttt{end}}=(L_{\texttt{end}},R_{\texttt{end}})$ denotes a final deck configuration, characterised by the fact that exactly one among $L_{\texttt{end}}$ and $R_{\texttt{end}}$ is empty.

To accommodate non-terminating matches in the formalism, without necessarily assuming that they exist, $\mathcal{S}^*$ is extended to
\begin{equation}
 \mathcal{S}^\infty\,:=\,\mathcal{S}^*\cup \mathcal{S}^\omega\,,
\end{equation}
where $\mathcal{S}^\omega$ contains infinite sequences.

 The trick-by-trick dynamics is governed by a \emph{trick rule function}
\begin{equation}
 \mathfrak{F}: \mathcal{S} \times \{A,B,\texttt{end}\} \to \mathcal{S} \times \{A,B,\texttt{end}\}
\end{equation}
that acts as follows:
\begin{itemize}
  \item For a playable start-of-trick configuration $(s,p)$, namely with $s = (L,R)$, $|L| \geqslant 1$, $|R| \geqslant 1$, $|L|+|R|=N$, and with $p \in \{A,B\}$ labelling the player moving first in the trick: $\mathfrak{F}$ executes one trick according to the game rules. If the resulting configuration has one empty deck (denoted as $s_{\texttt{end}}$, as already indicated above), then the output is $(s_{\texttt{end}},\texttt{end})$; otherwise the output is $(s',p')$ where $p'$ is the trick winner.
  \item $\mathfrak{F}(s,\texttt{end}) = (s,\texttt{end})$ (identity on terminal states).
  \item $\mathfrak{F}(s_{\texttt{end}},p) = (s_{\texttt{end}},\texttt{end})$ (conversion to terminal state).
\end{itemize}

 A match therefore evolves from $(s_0,A)$ along a deterministic trajectory of the type
 \begin{equation}
  (s_0,A)\:\mapsto\:\mathfrak{F}(s_0,A)\:\mapsto\:\mathfrak{F}^2(s_0,A)\:\mapsto\:\mathfrak{F}^3(s_0,A)\:\mapsto\:\cdots
 \end{equation}
 until possible game termination or, alternatively, entering a loop. In terms of the first-component projection
 \begin{equation}
  \mathcal{P}:\mathcal{S} \times \{A,B,\texttt{end}\}\to \mathcal{S}\,,\qquad \mathcal{P}(s,p)\,:=\,s\,,
 \end{equation}
 the corresponding finite or infinite sequence of decks configuration for that match is
 \begin{equation}
  (s_0,\mathcal{P}\mathfrak{F}(s_0,A),\mathcal{P}\mathfrak{F}^2(s_0,A),\dots)\,\in\,\mathcal{S}^\infty\,.
 \end{equation}

 Let $\mathcal{S}_{\textrm{start}}\subset\mathcal{S}$ denote the space of all evenly-split deck configurations of the $(N,\mathfrak{R})$-game, namely the start-of-match decks, i.e.,
 \begin{equation}
  \mathcal{S}_{\textrm{start}}\,:=\,\{s=(L,R)\in\mathcal{S}\,\textrm{ such that } |L|=|R|=N/2\}\,.
 \end{equation}
 Its (finite) cardinality, as counted already, is
 \begin{equation}
  |\mathcal{S}_{\textrm{start}}|\,=\,\mathcal{C}(N,\mathfrak{R})\,,
 \end{equation}
 and under the game rules (uniform shuffling followed by even split), each such arrangement has equal probability $1/\mathcal{C}(N,\mathfrak{R})$ to be selected at the beginning of a match.

 The match duration is naturally defined as the random variable
 \begin{equation}
  \begin{split}
     \mathcal{T}:\mathcal{S}_{\textrm{start}} & \to \mathbb{N} \cup \{\infty\}\,, \\
     \mathcal{T}(s)&:=\text{number of tricks in the match starting from configuration } s\,,
  \end{split}
 \end{equation}
  and let $\rho$ denote the probability distribution for $\mathcal{T}$, that is,
  \begin{equation}
   \begin{split}
    \rho(n) \,&:=\, \textrm{Prob}(\mathcal{T} = n) \,=\, \frac{|\{s \in \mathcal{S}_{\textrm{start}} : \mathcal{T}(s) = n\}|}{|\mathcal{S}_{\textrm{start}}|} \\
    &=\,\frac{\text{\,number of initial configurations yielding $n$-tricks-long matches}\,}{\mathcal{C}(N,\mathfrak{R})}\,.
   \end{split}
\end{equation}

Figures \ref{fig:stats_output_figure} and \ref{fig:r_parameter_plot} display numerical approximations of $\rho(n)$ empirically estimated from simulations over a (relatively small, in fact) subset of the total space of initial deck configurations. Since all simulated matches in that study terminated in finite time (i.e., with a finite number of tricks), these figures represent approximations of $\rho(n)$ for finite values $n \in \mathbb{N}$ only. The probability mass $\rho(\infty)$ corresponding to infinite-duration matches, if any exist, remains unobserved in those simulations.

 \section{Expectation of the game duration}\label{sec:expectation}

 In the notation of Section \ref{sec:abstr}, the average number of tricks per match is expressed by the expectation
\begin{equation}\label{eq:expectationTtricks}
 \mathbb{E}(\mathcal{T}) = \sum_{n\in\mathbb{N}} \,n \cdot \rho(n) + \infty \cdot \rho(\infty)
\end{equation}
 (with the tacit convention $\infty\cdot 0 = 0$). Actually, the series in \eqref{eq:expectationTtricks} is a finite sum. This follows from the fact that the number of distinct deck configurations is finite, $\mathcal{C}(N,\mathfrak{R})$, and consequently the set of possible durations for terminating matches is also finite. Thus, for fixed $N$ and $\mathfrak{R}$, there exists a maximum duration $\mathfrak{n}_{\max} \equiv \mathfrak{n}_{\max}(N,\mathfrak{R})$ such that no terminating match starting from any possible $(N,\mathfrak{R})$-deck configuration lasts longer than $\mathfrak{n}_{\max}$ tricks. (For concreteness, $\mathfrak{n}_{\max}$ can be chosen as the smallest integer with such property.) Thus, \eqref{eq:expectationTtricks} is re-written as
\begin{equation}\label{eq:expectationTtricks2}
\mathbb{E}(\mathcal{T}) = \sum_{n=1}^{\mathfrak{n}_{\max}} n \cdot \rho(n) + \infty \cdot \rho(\infty).
\end{equation}

Clearly, the numerical simulations reported in Figures \ref{fig:stats_output_figure}-\ref{fig:semilog_40_3} and Tables \ref{tab:numanalrecap} and \ref{tab:numanalrecapN}, only allow estimating the finite support sum  $\sum_{n=1}^{\mathfrak{n}_{\max}} n \cdot \rho(n)$ by using the empirical values for $\rho(n)$ and by cutting the sum to the maximum empirically observed match duration. This yields the (obviously finite) values for average number of tricks per match which are indicated in the tables above.

 From a theoretical perspective, the question of the finiteness of the game duration expectation was addressed first in \cite{Lakshtanov-Aleksenko-Beggar2013}. There, the finiteness of the mathematical expectation of the game duration is claimed to be proved, however under three fundamental departures from the game's classical rules:
 \begin{enumerate}
  \item each new trick is assumed to begin with a Bernoulli trial to select the starter player;
  \item at the end of each trick cards won by a player are assumed to be randomly shuffled before being added to the winner's deck;
  \item a further non-degeneracy condition is assumed, stating that given two players' decks at the beginning of a trick, the resulting updated decks at the end of the trick always differ depending on whether one or the other player started first.
 \end{enumerate}

 This creates a rather different game structure with no deterministic relationship between consecutive tricks and card ordering information destroyed at every trick, thus making the game a memory-less stochastic process rather than a deterministic sequence.

 Furthermore, whereas assumptions (1) and (2) above have theoretical interest per se, making the analysis of such `\emph{memory-less Beggar-My-Neighbour}' worthy, it is unclear whether assumption (3) is indeed applicable and holds for all card distribution: there could exist configurations where the trick outcome is the same regardless of starting player, an occurrence which is asserted in \cite{Lakshtanov-Aleksenko-Beggar2013} to be impossible, albeit without rigorous justification. In fact, assumption (3) is crucial for the graph-theoretic proof proposed in \cite{Lakshtanov-Aleksenko-Beggar2013}, for only by assuming that each state has exactly 2 distinct outgoing edges (corresponding to which player plays first) does the graph structure of the proof not collapse.

 Besides, we are unfortunately unable to reproduce a valid proof from the reasoning in \cite{Lakshtanov-Aleksenko-Beggar2013}, due to significant logical gaps and circular arguments.

 Despite the inconclusive nature of the approach in \cite{Lakshtanov-Aleksenko-Beggar2013}, what remains as a fundamental insight for our analysis is the direct connection between infinite matches and expectation finiteness.

 Indeed, following \eqref{eq:expectationTtricks2}, $\mathbb{E}(\mathcal{T})$ would be infinite if and only if any non-terminating match exists.
 If even one initial deck configuration produced a non-terminating match, that configuration having probability $1/\mathcal{C}(N,\mathfrak{R})>0$ under uniform shuffling, one would have $\rho(\infty)\geqslant 1/\mathcal{C}(N,\mathfrak{R})>0$, and consequently $\mathbb{E}(\mathcal{T}) =\infty$.

 The expectation $\mathbb{E}(\mathcal{T})$ is either finite (if no infinite matches exist) or infinite (if any infinite matches exist), with no intermediate possibility. The previous empirical findings of finite sums $\sum_{n=1}^{\mathfrak{n}_{\max}} n \cdot \rho(n)$ are therefore consistent with either scenario, since these sums exclude the infinite-duration contribution by construction.

 Clearly, the idea in \cite{Lakshtanov-Aleksenko-Beggar2013} underlying the stochastic version of the game is that randomisation eventually breaks cycles, and the `shuffling reset' at the end of each trick ensures eventual `absorption' (the graph of the match terminates). Instead, analysing the large-$k$ long-term behaviour of $\mathfrak{F}^k(s_0,A)$ requires understanding how card redistributions affect future accessibility, a combinatorially explosive problem.

\section{Forward vs backwards determinism}\label{sec:fwd-bkw-determinism}

The trick function $\mathfrak{F}$ is deterministic (in the customary sense, namely with only one possible output from each input), but it is not injective over $\mathcal{S}\times\{A,B\}$. That is, there are distinct configurations for the start of a trick, say, $(s,p)\neq(s',p')$, which yield the same configuration $\mathfrak{F}(s,p)=\mathfrak{F}(s',p')$ at the end of that trick.

This fact appears to have been highlighted explicitly for the first time in the recent work \cite{Casella2024}. Here we present a more systematic perspective.

The lack of injectivity of $\mathfrak{F}$ manifests at various levels.
\begin{enumerate}
 \item[1.] $\mathfrak{F}\upharpoonright  (\mathcal{S}\times\{A,B\})_{\textrm{playable}}$ is not injective. That is, $\mathfrak{F}$ is not injective even when restricted to the collection of start-of-trick configurations
\begin{equation}
 (\mathcal{S}\times\{A,B\})_{\textrm{playable}}\,:=\,\big\{(s,p)\in \mathcal{S}\times\{A,B\}\textrm{ such that } s=(L,R) \textrm{ with }|L|,|R|\geqslant 1\big\}\,.
\end{equation}
\item[2.] $\mathfrak{F}\upharpoonright  (\mathcal{S}\times\{A,B\})_{\textrm{reachable}}$ is not injective. That is, $\mathfrak{F}$ is not injective even when further restricted to the sub-class of start-of-trick configurations that can be actually \emph{reached} from one of the $\mathcal{C}(N,\mathfrak{R})$ evenly split initial decks, namely
  \begin{equation}\label{eq:reachablespace}
 \begin{split}
   & (\mathcal{S}\times\{A,B\})_{\textrm{reachable}}\,:= \\
   &\quad :=\,\left\{(s,p)\in (\mathcal{S}\times\{A,B\})_{\textrm{playable}}\,\textrm{ such that }
 \begin{array}{c}
   (s,p)=\mathfrak{F}^m(s_0,A)\textrm{ for some }m\in\mathbb{N}_0, \\
   \textrm{and some initial-match deck $s_0$}
 \end{array}\!
 \right\}.
 \end{split}
\end{equation}
\end{enumerate}
We observe furthermore that
  \begin{equation}\label{eq:strictInclTrickReach}
  (\mathcal{S}\times\{A,B\})_{\textrm{reachable}}\,\varsubsetneq\,(\mathcal{S}\times\{A,B\})_{\textrm{playable}}\,.
 \end{equation}

 Verifying the first claim is straightforward. Among other similar examples, consider the following mechanism: in either starting configuration
\begin{equation}\label{eq:nonInj1}
\begin{split}
\texttt{DeckA} &= \texttt{[1CC}\boxed{\texttt{*\dots*}}\texttt{C2]}\,, \\
\texttt{DeckB} &= \texttt{[C1C]}
\end{split}
\end{equation}
and
\begin{equation}\label{eq:nonInj2}
\begin{split}
\texttt{DeckA} &= \texttt{[C1CC}\boxed{\texttt{*\dots*}}\texttt{]}\,, \\
\texttt{DeckB} &= \texttt{[2C1C]}\,,
\end{split}
\end{equation}
player $A$ moves first and wins the trick, outputting the decks
\begin{equation}\label{eq:nonInj3}
\begin{split}
\texttt{DeckA} &= \texttt{[CC}\boxed{\texttt{*\dots*}}\texttt{C21C]}\,, \\
\texttt{DeckB} &= \texttt{[1C]}\,.
\end{split}
\end{equation}
Here, $\boxed{\texttt{*\dots*}}$ is the same sequence of $N-8$ cards: the argument is applicable to decks of any size and can be naturally iterated to produce multiple distinct start-of-trick decks with the same end-of-trick output.

\begin{figure}[!htbp]
   \centering
   \includegraphics[width=1\textwidth]{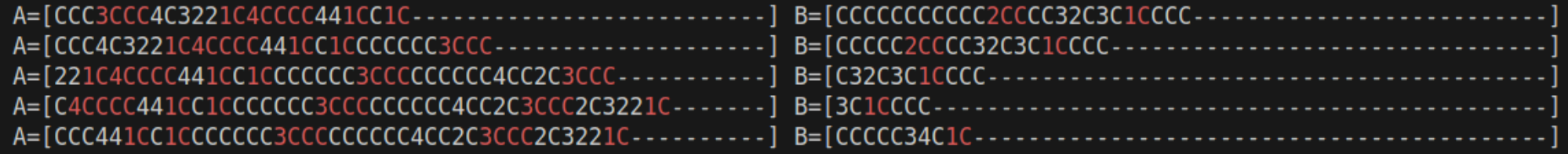}

   \bigskip

   \includegraphics[width=1\textwidth]{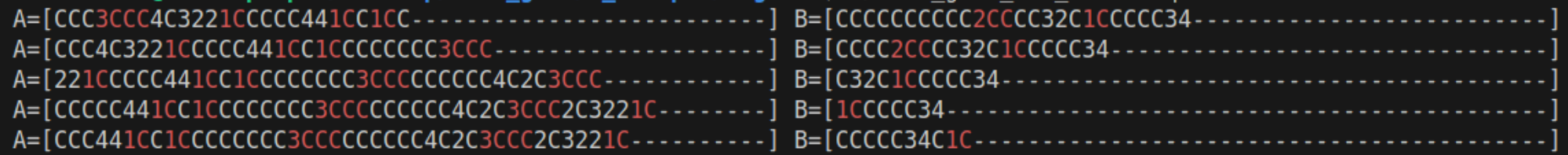}
   \caption{\small Explicit detection of lack of injectivity, in the space of reachable configurations, of the trick rule function $\mathfrak{F}$ in the $(52,4)$-setting: two distinct matches display distinct decks after the third trick, and yet the same decks in the following fourth trick.  \normalsize}
   \label{fig:noninj}
\end{figure}

As for the second claim, counterexamples have been identified, albeit not systematically over all $(N,\mathfrak{R})$-scenarios, to the injectivity of $\mathfrak{F}\upharpoonright  (\mathcal{S}\times\{A,B\})_{\textrm{reachable}}$. For instance, elaborating from the analysis of \cite{Casella2024} in the $(52,4)$-setting, Figure \ref{fig:noninj} reports the first 4 tricks of the dynamics of two distinct matches originating from initial configurations $s_0$ and $\widetilde{s}_0$, $s_0\neq\widetilde{s}_0$; thus all subsequent $\mathfrak{F}^n(s_0,A)$'s and $\mathfrak{F}^n(\widetilde{s}_0,A)$'s, for $n\in\{0,1,2,3,4\}$, surely belong to $(\mathcal{S}\times\{A,B\})_{\textrm{reachable}}$. One sees by inspection that
\begin{equation}
 \mathfrak{F}^3(s_0,A)\neq\mathfrak{F}^3(\widetilde{s}_0,A)\,\quad\text{ and yet }\quad \mathfrak{F}^4(s_0,A)=\mathfrak{F}^4(\widetilde{s}_0,A)\,,
\end{equation}
which shows lack of injectivity of $\mathfrak{F}$ in $(\mathcal{S}\times\{A,B\})_{\textrm{reachable}}$.

As a consequence of $\mathfrak{F}\upharpoonright  (\mathcal{S}\times\{A,B\})_{\textrm{reachable}}$ failing to be injective, there is no global backwards determinism in the space of reachable configurations: reconstructing the decks that yielded a given trick's output produces in general a backward analysis of the match through an expanding tree of former decks, rather than a unique history.

    Finally, let us establish the strict inclusion \eqref{eq:strictInclTrickReach}: not all start-of-trick configurations can be part of an actual match. A necessary condition for a still playable configuration $(s,p)$ to be a trick's output in an $(N,\mathfrak{R})$-game, starting from a playable trick's input, is that one of the two players' decks in $s=(L,R)$, and precisely player $p$'s deck (the trick's winner), ends with a block of the type
\begin{equation}\label{eq:endingblocks}
 \texttt{[}k \underbrace{\texttt{C\dots C}}_{(k \text{ times})}\texttt{]}\,,\qquad k\in\{1,\dots,\mathfrak{R}\}\,,
\end{equation}
that is, for concreteness, bottom blocks like
\[
 \texttt{[1C]}\,,\qquad \texttt{[2CC]}\,,\qquad \texttt{[3CCC]}\,,\qquad \textrm{etc.}
\]
Indeed, player $p$'s victory in a trick can eventually only occur by playing a rank-$k$ card to which the other player is unable to respond in $k$ moves. Figure \ref{fig:noninj} highlights this fact in red. It is easy to violate condition \eqref{eq:endingblocks} and generate deck configurations of the type
\begin{equation}
\begin{split}
\texttt{DeckA} &= \texttt{[*\dots*}\:k\underbrace{\texttt{C\dots C}}_{(\neq k \text{ times})}\texttt{]}\,, \\
\texttt{DeckB} &= \texttt{[*\dots*}\:k'\underbrace{\texttt{C\dots C}}_{(\neq k' \text{ times})}\texttt{]}\,,
\end{split}
\end{equation}
for $k,k'\in\{1,\dots,\mathfrak{R}\}$, which give start-of-trick configurations with suitable choice of the variable part \texttt{[*, $\cdots$, *]}, and yet cannot have been generated at the end of a trick and therefore cannot be part of an actual match.

Observe that an empirical indication of \eqref{eq:strictInclTrickReach} could have been obtained by considering the cardinalities
 \begin{equation}
 \begin{split}
  \big| \mathcal{S}\times\{A,B\}\big|\,&=\,2(N+1)\,\mathcal{C}(N,\mathfrak{R})\,, \\
  \big| (\mathcal{S}\times\{A,B\})_{\textrm{playable}} \big|\,&=\,2(N-1)\,\mathcal{C}(N,\mathfrak{R})\,.
 \end{split}
\end{equation}
They are determined by multiplying the number $\mathcal{C}(N,\mathfrak{R})$ of distinct card permutations (equivalent to the count of distinct initial decks) by the number of valid splits for each permutation (respectively $N+1$ for the total space and $N-1$ for the playable space, and further by 2 accounting for the two players).
Then, at least for those simulations reported in Tables \ref{tab:numanalrecap} and \ref{tab:numanalrecapN} for which the empirical average $\overline{\tau}$ of the number of tricks per match is considerably smaller than $2(N-1)$, an (over)estimate of the number of reachable tricks is possible by counting each match's tricks as different from the other matches', which yields
 \begin{equation}\label{eq:estimateCard}
   \big| (\mathcal{S}\times\{A,B\})_{\textrm{reachable}} \big|\,\lesssim\,\overline{\tau}\cdot\mathcal{C}(N,\mathfrak{R})\,<\,\big| (\mathcal{S}\times\{A,B\})_{\textrm{playable}} \big|\,.
 \end{equation}
 For example \eqref{eq:estimateCard} is applicable to the scenario $N=40$,  since $\overline{\tau}\leqslant 41.26$ versus $2(N-1)=78$ (Table~\ref{tab:numanalrecap}), as well as to the scenario $N=52$, since $\overline{\tau}=51.16$ versus $2(N-1)=102$ (Table~\ref{tab:numanalrecapN}).

\section{Patterns and distributions of ultra-long matches}\label{sec:et}

  There appears to have been a widespread activity, from the amateur level onward, in identifying initial deck configurations that give rise to ultra-long matches of Beggar-My-Neighbour and its variants, such as \cite{Paulhus_Beggar1999,Beggar-Collins-records,Beggar-Mann-records,Beggar-Mayer-codes-records,Beggar-Zanotto-codes,Beggar-Tristan-records,Beggar-Gentilini}, or the attempts reported in \cite{Casella2024}, and many more. By ultra-long one may conventionally mean matches whose length exceeds the average length by many standard deviations. Clearly, in this context by average length we refer to the component $\sum_{n=1}^{\mathfrak{n}_{\max}} n \cdot \rho(n)$ of $\mathbb{E}(\mathcal{T})$, in the notation of \eqref{eq:expectationTtricks2}.

  This activity has primarily focused on the search for progressively longer, terminating matches. However, a systematic analysis of the structure and patterns within these ultra-long matches is currently lacking. This Section aims to outline such an analysis.

  Furthermore, whereas we could reproduce with our simulation scripts the ultra-long matches reported, for instance, in \cite{Beggar-Collins-records,Beggar-Mann-records,Beggar-Mayer-codes-records,Beggar-Tristan-records}, there are several other amateur sources indicating initial decks that actually fail to yield the declared ultra-long matches. For certain authors with accessible scripts (as \cite{Beggar-Gentilini}) we could ascertain that such a discrepancy is due to the numerical implementation of a variant of the game, where at the end of each trick, when the challenge is concluded, the pile of cards sitting face-up on table are picked one by one from top to bottom, and each goes face-down to bottom of winner's deck, with the net result that the card order from that trick is reversed. In the following we shall discard this type of minority variant. In other cases where the used codes are inaccessible, we simply cannot reproduce the announced results.

  For the subsequent analysis of ultra-long matches, we therefore relied on the well-validated numerical results presented in Tables \ref{tab:numanalrecap} and \ref{tab:numanalrecapN}, along with other ultra-long matches that we were able to reproduce and verify from external sources, all concerning various $(N,\mathfrak{R})$-settings. Noticeably, for smaller $\mathfrak{R}/N$-ratios both the average number of tricks per match and their dispersion dramatically increase, resulting in increasingly longer matches.

 For each considered ultra-long match, the match history was systematically screened. After each completed trick, the following parameters were recorded for each player's deck:
 \begin{itemize}
  \item total number of cards in the deck;
  \item number of special cards contained within the deck;
  \item the positions (counted from the top) of all special cards.
 \end{itemize}


 \bigskip

\begin{figure}[!ht]
    \centering

    \includegraphics[width=1\textwidth]{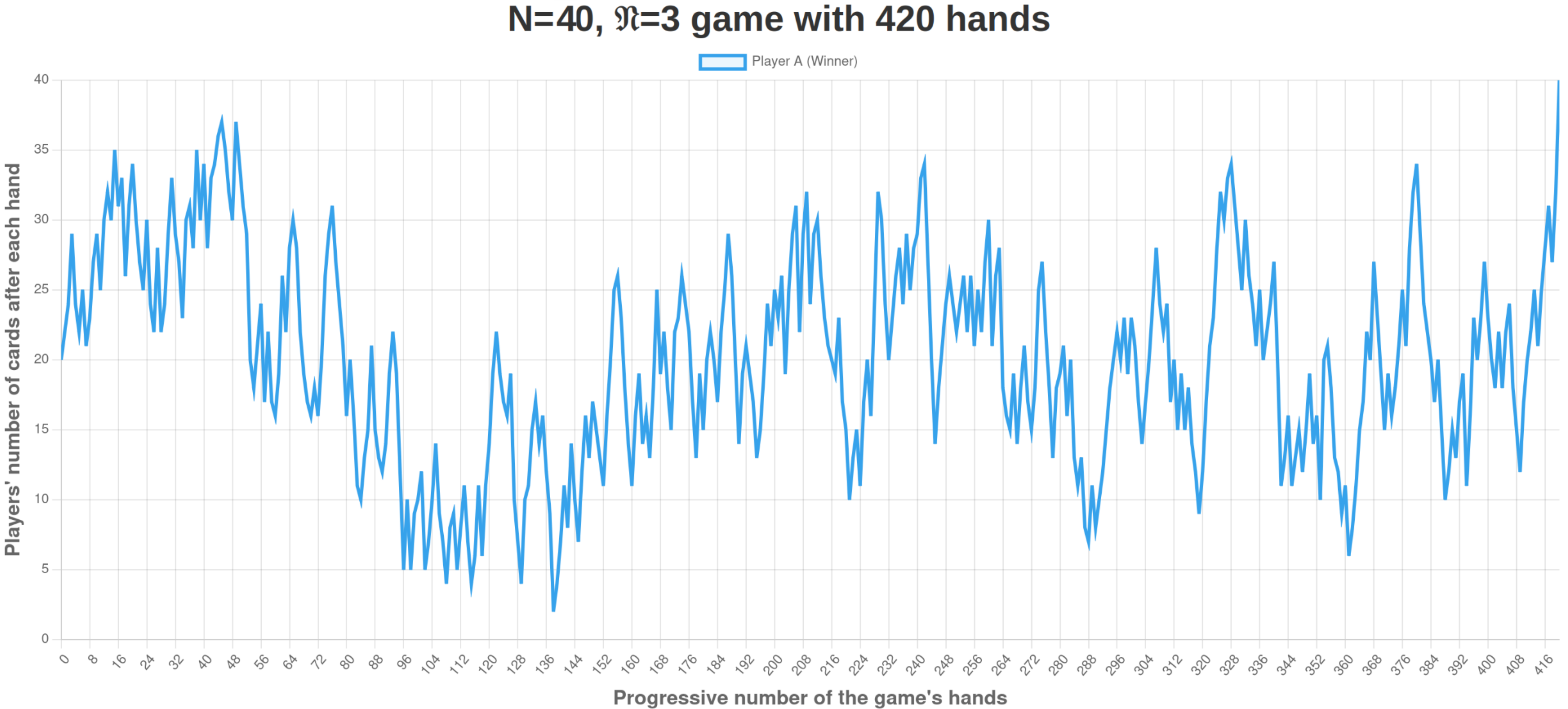}

     \bigskip

    \includegraphics[width=1\textwidth]{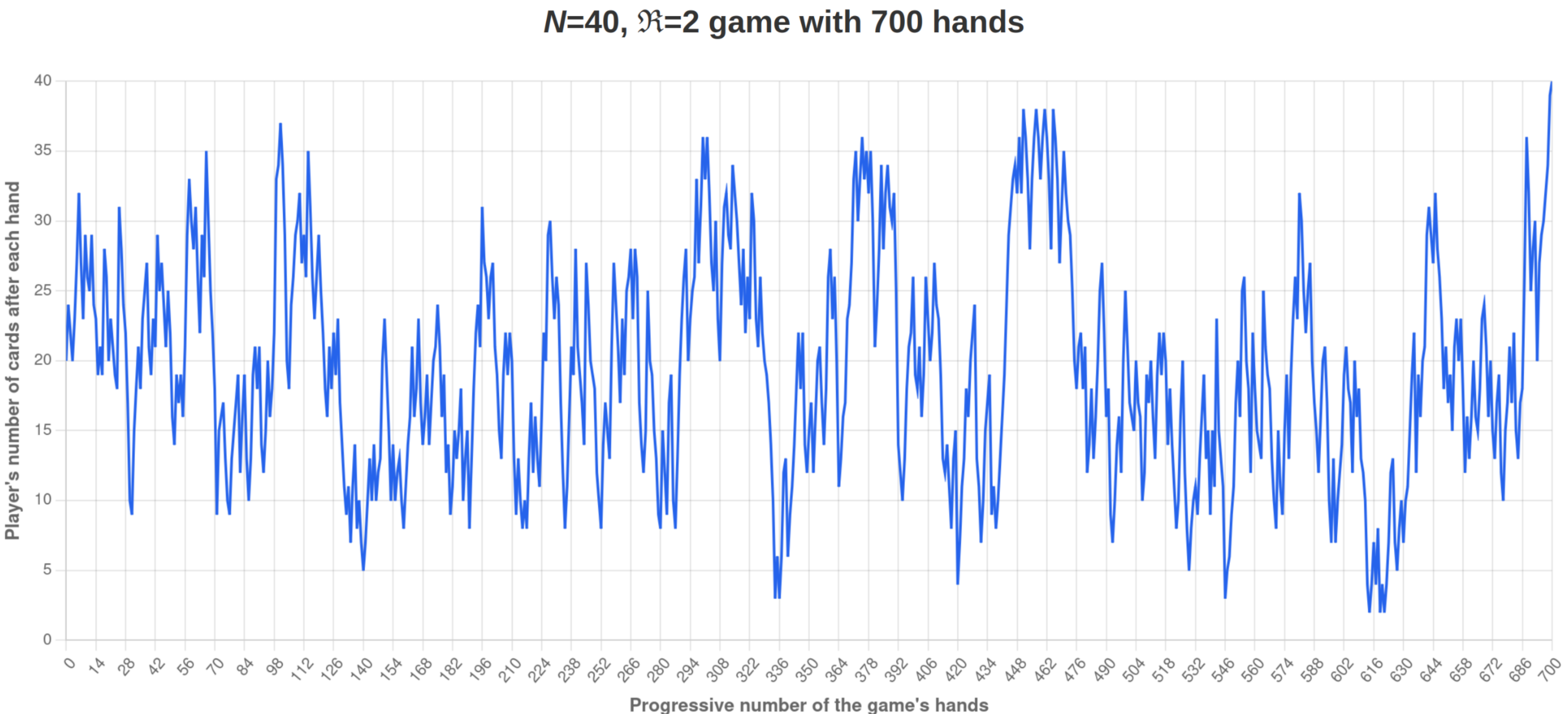}

    \caption{\small Dynamics of two ultra-long matches detected in the numerical simulation. Top: a 420 tricks long $(40,3)$-game won by player $A$ with starting decks \eqref{eq:startdeck-40-3-420}. Bottom: a 700 tricks long $(40,2)$-game won by player $B$ with starting decks \eqref{eq:startdeck-40-2-700}. \normalsize}
    \label{fig:tricks_oscillation_plot_R3}
\end{figure}

 $\bullet$ \textbf{Pseudo-periodic patterns in the number of cards.} From the simulations summarised in Tables \ref{tab:numanalrecap} and \ref{tab:numanalrecapN}, Figure \ref{fig:tricks_oscillation_plot_R3} visualises the dynamics of
\begin{itemize}
 \item a 420-tricks long $(40,3)$-game won by player $A$ with starting decks
\begin{equation}\label{eq:startdeck-40-3-420}
\begin{split}
\texttt{DeckA} &= \texttt{[C1CC1CC32C3CCC3C2CCC]}\,, \\
\texttt{DeckB} &= \texttt{[CC2CCC2CC1CC3CC1CCCC]}
\end{split}
\end{equation}
(this is actually the longest match previously reported in Figure \ref{fig:stats_output_figure}),
 \item and a 700-tricks long $(40,2)$-game won by player $B$ with starting decks
 \begin{equation}\label{eq:startdeck-40-2-700}
\begin{split}
\texttt{DeckA} &= \texttt{[CCCC1C1CCCCCCCCCCCCC]}\,, \\
\texttt{DeckB} &= \texttt{[C22CCCC1CC2CCCC1CCC2]}\,.
\end{split}
\end{equation}
\end{itemize}

\begin{figure}[!ht]
    \centering

    \includegraphics[width=1\textwidth]{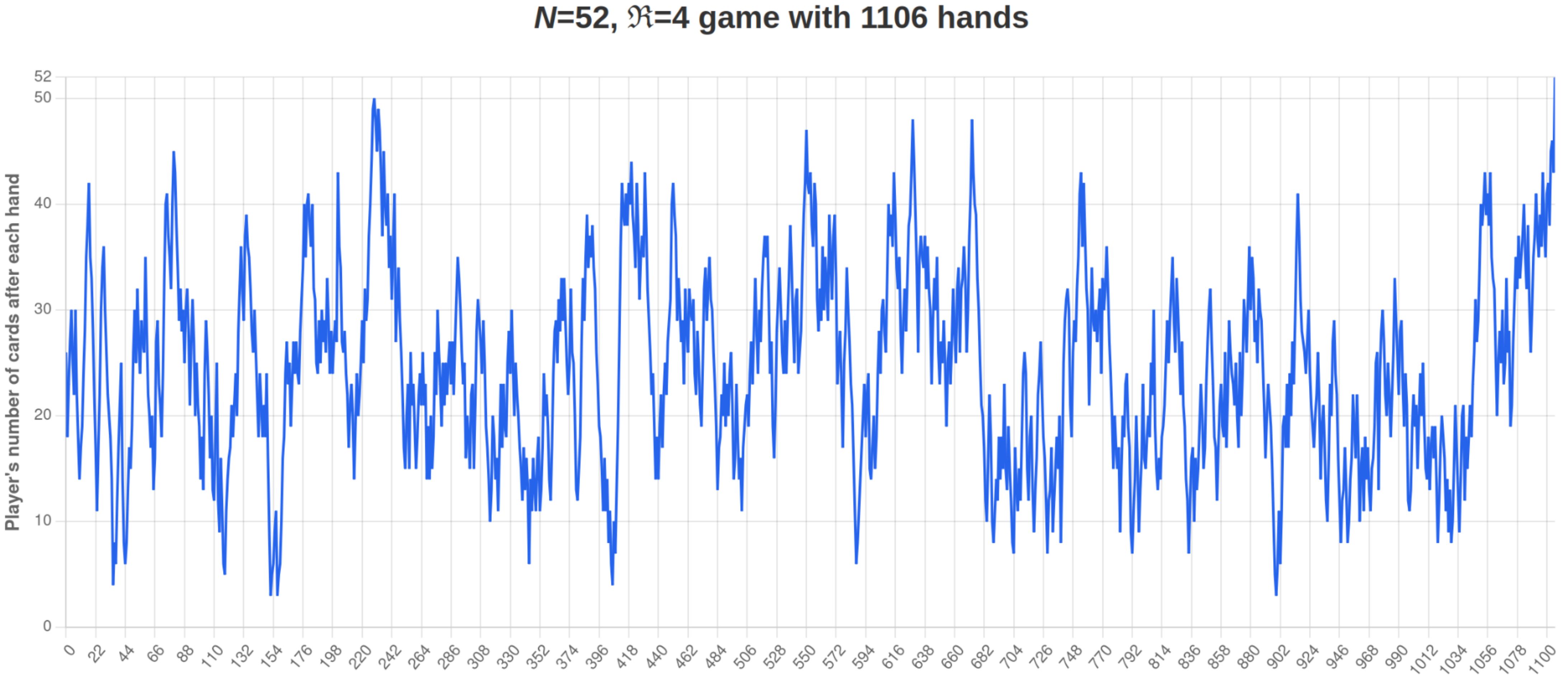}

     \bigskip

    \includegraphics[width=1\textwidth]{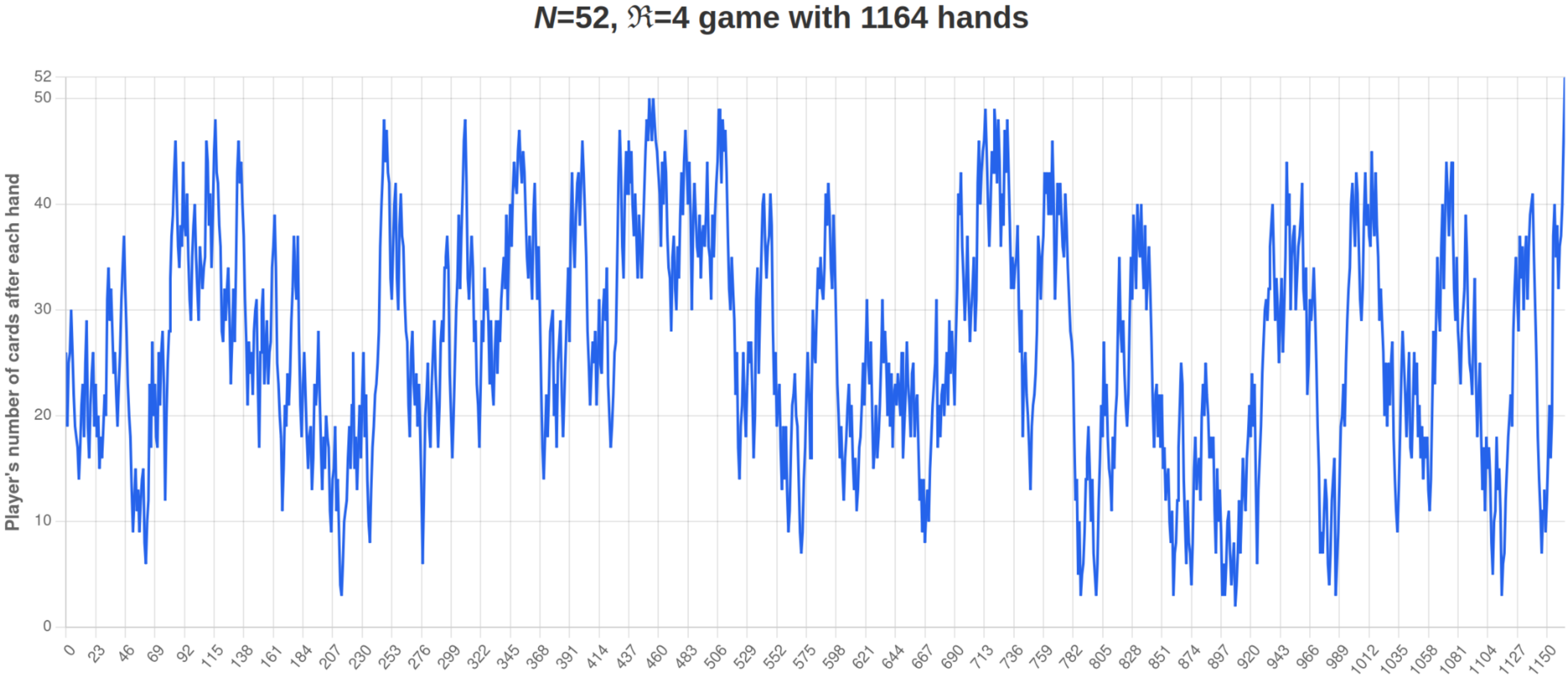}

    \caption{\small Dynamics of two ultra-long matches detected in the numerical simulation. Top: a 1106 tricks long $(52,4)$-game won by player $A$ with starting decks \eqref{eq:startdeck-52-4-1106}. Bottom: a 1164 tricks long $(52,4)$-game won by player $B$ with starting decks \eqref{eq:startdeck-52-4-1164}. \normalsize}
    \label{fig:tricks_oscillation_plot_52-4}
\end{figure}

 Analogously, Figure \ref{fig:tricks_oscillation_plot_52-4} visualises the dynamics of
\begin{itemize}
 \item a 1106-tricks long $(52,4)$-game won by player $A$ with starting decks
\begin{equation}\label{eq:startdeck-52-4-1106}
\begin{split}
\texttt{DeckA} &= \texttt{[CCCC3CCC4CC2C4CC114CCCCCC1]}\,, \\
\texttt{DeckB} &= \texttt{[CCCCC33CCCCCCCCC4C13C2C2C2]}\,,
\end{split}
\end{equation}
 \item and a 1164-tricks long $(52,4)$-game won by player $B$ with starting decks
 \begin{equation}\label{eq:startdeck-52-4-1164}
\begin{split}
\texttt{DeckA} &= \texttt{[CCC41CC2CCCCCCCCC243211C23]}\,, \\
\texttt{DeckB} &= \texttt{[CCCCC4CCCC31C3CCCCCCCC4CCC]}\,,
\end{split}
\end{equation}
\end{itemize}
 whose identification was reported, among others, in \cite{Beggar-Mann-records,Beggar-Mayer-codes-records,Casella2024}.

 All such matches exhibit a consistent multi-scale oscillatory structure characterised by $(N,\mathfrak{R})$-dependent temporal hierarchies. Specifically:
\begin{itemize}
\item \emph{Macro-cycles} showing major dominance shifts between players, with characteristic period increasing dramatically as $\mathfrak{R}/N$ decreases: approximately $80$-$100$ tricks for the $(40,3)$-game, $120$-$150$ tricks for the $(40,2)$-game, while the $(52,4)$-games exhibit more frequent but smaller-amplitude oscillations with periods of approximately $7$ tricks between successive peaks.

\item \emph{Meso-cycles} of intermediate-scale fluctuations of $20$-$30$ tricks, representing secondary oscillatory dynamics nested within the macro-cycles, particularly evident in the lower $\mathfrak{R}/N$ matches.

\item \emph{Micro-oscillations} of high-frequency variations of $\pm 2$ to $\pm 5$ cards, reflecting individual trick-to-trick dynamics, and present across all matches regardless of parameters.
\end{itemize}

  \bigskip

 $\bullet$ \textbf{Approximate geometric regime.} To a partial extent, such behaviour in the somewhat erratic, stock-market-like dynamics of the decks' cardinalities along ultra-long matches is understandable in terms of the approximate exponential distribution of match duration (Section \ref{sec:num}, and Figure \ref{fig:semilog_40_3}). This line of reasoning is already present in \cite{Casella2024} in a form that deserves the following clarifications.

 The numerical evidence that $\rho(n)\approx e^{-\lambda n}$ in wide intermediate regimes of $n$ suggests considering an idealised situation (realistic in those regimes) where the distribution $\rho$ is precisely exponential throughout its domain $\mathbb{N}$ (having temporarily factored out the expectedly super rare occurrence $n=\infty$), namely
 \[
  \rho(n)\,=\, C e^{-\lambda n}\,,\quad n\in\mathbb{N}\,.
 \]
 It is standard to argue that the normalisation condition $\sum_{n\in\mathbb{N}}\rho(n)=1$ determines $C=\frac{1-e^{-\lambda}}{e^{-\lambda}}$, and setting $p:=1-e^{-\lambda}$ yields
 \[
  \rho(n)\,=\, \frac{1-e^{-\lambda}}{e^{-\lambda}} \cdot e^{-\lambda n} = (1-e^{-\lambda}) e^{-\lambda(n-1)}\,=\,p(1-p)^{n-1}\,,\quad n\in\mathbb{N}\,.
 \]

  This is the geometric distribution with success probability $p$ (the discrete counterpart of the continuous exponential distribution). As a geometric random variable, the match duration $\mathcal{T}$ displays well-known memory-less properties, say,
  \begin{equation}\label{eq:memory-less}
   \textrm{Prob}(\mathcal{T} > m+n \mid \mathcal{T} > m) \,=\, \textrm{Prob}(\mathcal{T} > n)\quad \forall m,n\in\mathbb{N}\,,
  \end{equation}
 and in particular the conditional probability of match termination at the next trick is constant,
 \begin{equation}\label{eq:constant_collapse}
\textrm{Prob}(\mathcal{T} = n+1 \mid \mathcal{T} > n) = p\quad \forall n\in\mathbb{N}\,.
\end{equation}

 Here \eqref{eq:constant_collapse} refers to sampling from the ensemble $\mathcal{S}_{\textrm{start}}$ of initial configurations (sampling $s$ uniformly from $\mathcal{S}_{\textrm{start}}$ conditional on the event $\{\mathcal{T}(s) > n\}$, i.e., from configurations yielding games lasting at least $n$ tricks, then the probability that $\mathcal{T}(\omega) = n+1$ is $p$, independent of $n$). It is \emph{not} a statement about stochasticity within a single match (which is, instead, fully deterministic). So the statement of \cite{Casella2024} that ``any ongoing game has a constant probability of `collapsing' as it progresses'' must be carefully placed into this context: it is only indirectly justified by the overall probabilistic description of the game, and does not describe any random mechanism within the match dynamics. The memory-less property \eqref{eq:memory-less} should therefore be understood as an approximation valid for the typical duration range, not as an exact property across all game lengths.

 A more accurate, related comment from \cite{Casella2024} concerns the non-predictive power of past durations: ``the fact that a game has endured for a certain amount of time has no predictive power over its future length''. Indeed, the memory-less property \eqref{eq:memory-less} implies:
\begin{equation}\label{eq:memory-less3}
\mathbb{E}[\,\mathcal{T} - n \mid \mathcal{T} > n] = \mathbb{E}[\mathcal{T}] = \frac{1}{p}\,,
\end{equation}
 that is, the expected number of remaining tricks is always $1/p$, regardless of how many tricks have already occurred. Conditioning on survival to time $n$ does not change the distribution of remaining lifetime when sampling from the $\mathcal{S}_{\textrm{start}}$  ensemble.

 The geometric model is itself an approximation, as said. It would predict mean and variance as
 \begin{equation}
  \mathbb{E}(\mathcal{T})\,=\,\frac{1}{p}\,,\qquad \textrm{Var}(\mathcal{T})\,=\,\frac{1-p}{p^2}\,,
 \end{equation}
 hence variance-to-mean ratio
\begin{equation}
\frac{\text{Var}(\mathcal{T})}{\mathbb{E}(\mathcal{T})} \,=\, \frac{1-p}{p} \,=\, \mathbb{E}(\mathcal{T}) - 1 \,\approx\, \mathbb{E}(\mathcal{T})
\end{equation}
for small $p$ (long games). Our Tables \ref{tab:numanalrecap} and \ref{tab:numanalrecapN} show variance-to-mean ratios ranging from $0.89$ to $37.76$, which deviate significantly from this prediction: in all considered settings the observed ratio remains below the geometric benchmark $\mathbb{E}(\mathcal{T})-1$, i.e., the data are under-dispersed relative to the purely geometric model. Combined with the Poisson-baseline analysis of Section \ref{sec:num} (over-dispersion, ratio $>1$, for low $\mathfrak{R}/N$; under-dispersion, ratio $<1$, for high $\mathfrak{R}/N$), this indicates correlations, clustering effects, and regulatory constraints beyond the simple geometric model. The pure exponential/geometric behaviour may only emerge in intermediate regimes or as a coarse approximation.

The approximation of purely geometric game duration $\mathcal{T}$ with parameter $p = 1 - e^{-\lambda}$, where $\lambda \approx 1/\mathbb{E}(\mathcal{T})$, is also useful for an order-of-magnitude assessment of the probability of observing the documented ultra-long matches. For the $(40,3)$-setting with average duration $\mathbb{E}(\mathcal{T}) \approx 30.61$ tricks (Table \ref{tab:numanalrecap}), using the fitted $\lambda \approx 0.0394$ (Figure \ref{fig:semilog_40_3}) and hence $p \approx 0.0386$, the probability of a match lasting at least 420 tricks would be
\begin{equation}
\textrm{Prob}(\mathcal{T} \geqslant 420) \,= \,(1-p)^{419} \,\approx\, e^{-\lambda \cdot 419} \,\approx\, e^{-16.5} \,\approx\, 6.8\cdot 10^{-8}\,.
\end{equation}
  Similarly, for the $(40,2)$-setting with average $\mathbb{E}(\mathcal{T})\approx 41.26$ tricks, a 700-trick match would have probability:
 \begin{equation}
\textrm{Prob}(\mathcal{T} \geqslant 700) \,\approx\, e^{-700/41.26} \,\approx\, e^{-17.0} \,\approx\, 4.1\cdot10^{-8}\,.
\end{equation}
 So, with $10^7$ simulated matches, the expected number of games lasting 420 tricks or more would be $10^7 \cdot 6.8\cdot10^{-8}\approx 0.7$ games. And analogously, with $10^8$ simulated matches, the expected number of games lasting 700 tricks or more would be $10^8 \cdot 4.1\cdot10^{-8}\approx 4$ games.
Having observed such games in the documented simulations is consistent with the exponential approximation for $\rho(n)$, suggesting that ultra-long matches, while rare, do occur at rates captured by this simple model.

   \bigskip

 $\bullet$ \textbf{Further non-geometric patterns in the number of cards.} Further distinguished features of ultra-long matches are in common across different $(N,\mathfrak{R})$-settings.
\begin{itemize}
 \item The oscillatory behaviour remains remarkably stable throughout the vast majority of each match history, with the card count oscillating around the natural balance point $N/2$ (`\emph{sustained equilibrium phase}'). 
 \item There occur multiple \emph{`near-death recoveries'}, where one player's deck diminishes to critically low levels ($2$ cards in the $(40,3)$-game at trick $138$; $2$-$3$ cards multiple times in the $(40,2)$-game around tricks $450$-$470$; $3$ cards in the $(52,4)$ $1106$-tricks game), followed by remarkable recoveries.
 \item Recovery trajectories from extreme lows seem to follow gradual, irregular paths, whereas declines from peaks ($37$ cards at trick $49$ in the $(40,3)$-game; $39$ cards near trick $700$ in the $(40,2)$-game; peaks around $50$ in both $(52,4)$-games) exhibit sharper, more direct downward movements (`\emph{asymmetric mean reversion}').
 \item Matches eventually terminate with a `\emph{rapid terminal resolution}' in which the balance shifts decisively only near the end of each match, with victory occurring relatively swiftly after hundreds of tricks in quasi-equilibrium.
\end{itemize}

 All these are implicit signatures of deviation from the purely geometric/exponential approximation.


   \bigskip

 $\bullet$ \textbf{Patterns in the initial configuration effects.} A first indication is that ultra-long matches can arise from diverse initial special card distributions. The $(40,3)$-decks \eqref{eq:startdeck-40-3-420} show nearly balanced special cards between players in both number and rank. In contrast, the $(40,2)$-game with decks \eqref{eq:startdeck-40-2-700} exhibits significant asymmetry with player $B$ holding $6$ special cards versus player $A$'s $2$, yet player $A$ sustained the match for $700$ tricks before ultimately losing. The $(52,4)$-examples show intermediate distributions. This suggests that strategic positioning and spacing of special cards within the deck ordering can substantially compensate for numerical disadvantages in special card count. The deterministic game dynamics allow even initially disadvantaged configurations to generate complex, persistent oscillatory behaviour where eventual victory depends sensitively on the detailed card sequencing rather than merely on the gross inventory of special cards.

   \bigskip

 $\bullet$ \textbf{Patterns of the special cards separation.} The mean separation among special cards in the winner's final deck was computed over many long matches. This is in fact a combined measure of the numbers of all cards and of special cards. When a player's deck has size $m$ and contains $k$ special cards occupying positions $p_1, \ldots, p_k$ (from top to bottom, hence $1 \leqslant p_1 < p_2 < \cdots < p_k \leqslant m$), the special cards partition the deck into $k+1$ intervals
\begin{equation}\label{eq:intervalsSEP}
I_0 = [1, p_1), \quad I_j = (p_j, p_{j+1}), \; j \in \{1, \ldots, k-1\}, \quad I_k = (p_k, m],
\end{equation}
with lengths $\ell_0, \ell_1, \ldots, \ell_k$ (counting only ordinary cards), that is,
\begin{equation}\label{eq:ellsSEP}
\begin{split}
\ell_0 \,&=\, p_1 - 1\,, \\
\ell_j \,&=\, p_{j+1} - p_j - 1, \quad j \in \{1, \ldots, k-1\}\,, \\
\ell_k\, &=\, m - p_k\,.
\end{split}
\end{equation}
Thus, $\sum_{j=0}^{k} \ell_j = m - k$ (total ordinary cards), and
 \begin{equation}\label{eq:meansep}
  \textrm{mean separation}\, =\,\frac{1}{k+1}\sum_{j=0}^k \ell_j\,=\,\frac{m-k}{k+1}\,.
 \end{equation}

\begin{figure}[!htbp]
    \centering
    \moveleft1.8cm\vbox{%
        \hbox{\includegraphics[width=0.6\textwidth]{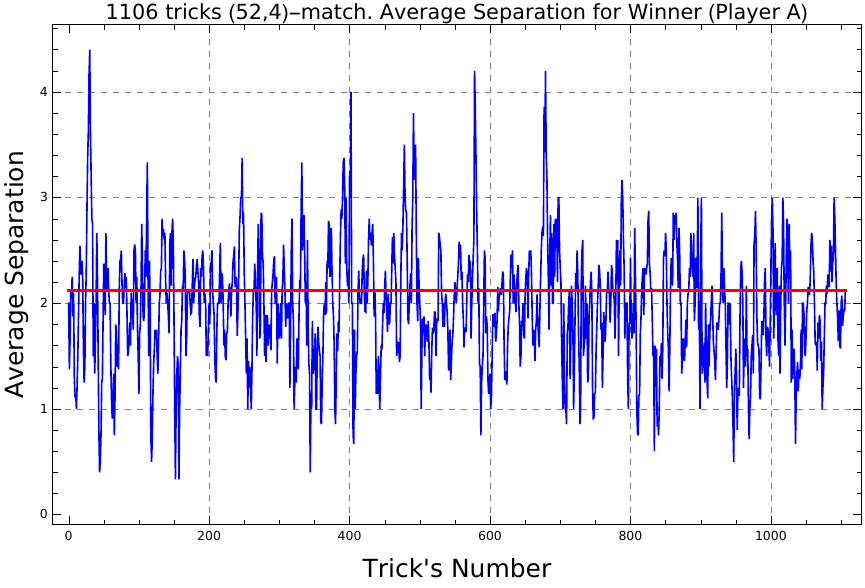}\quad\includegraphics[width=0.6\textwidth]{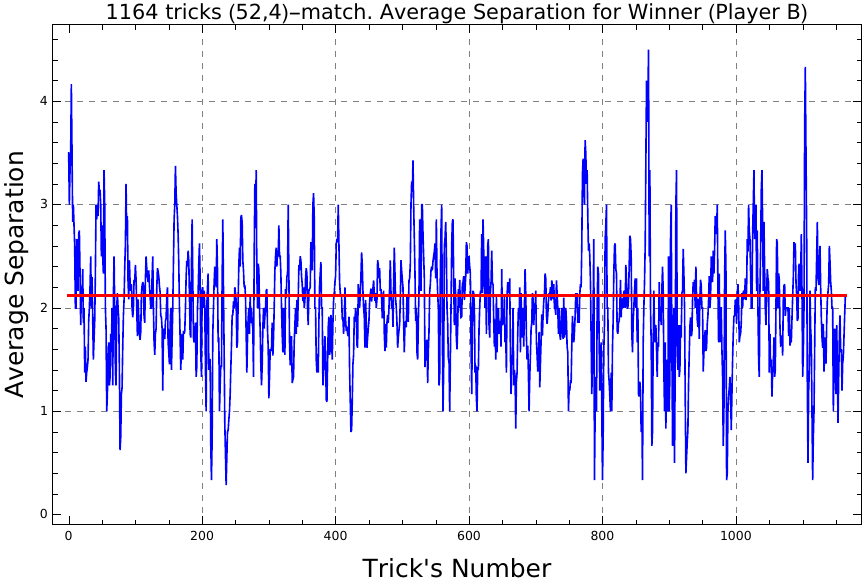}}%
    }
    \caption{\small Evolution of the average separation among special cards in the deck of the winner for two ultra-long $(52,4)$-matches. \normalsize}
    \label{fig:separation}
\end{figure}

  Numerics on ultra-long matches reveal that while card counts exhibit multi-scale oscillations that eventually resolve to victory, mean separation fluctuates around the equilibrium value $(N-4\mathfrak{R})/(4\mathfrak{R}+1)$ throughout the entire match (Figure~\ref{fig:separation}), reaching extremes pseudo-periodically.
 \begin{itemize}
  \item High separation indicates weak tactical position: few special cards ($k+1$ small in \eqref{eq:meansep}) dispersed among many ordinary cards ($m-k$ large).
  \item Low separation indicates strong tactical position: many special cards ($k$ large) densely packed with few ordinary cards between them.
 \end{itemize}
The oscillations exhibit amplitude asymmetry -- downwards excursions are more extreme than upwards spikes -- and maintain the same average amplitude and frequency throughout the match, a behaviour typical of random processes, even though the card redistribution is entirely deterministic. The separation converges exactly to $(N-4\mathfrak{R})/(4\mathfrak{R}+1)$ at match termination when the winner holds all cards.

   \bigskip

 $\bullet$ \textbf{Patterns in the position entropy.} Closely related to the special cards' average separation is the position entropy $H$ defined as
\begin{equation}\label{eq:defPosEntr}
H := -\sum_{\substack{j=0 \\ \ell_j > 0}}^{k} \frac{\ell_j}{m} \ln\left(\frac{\ell_j}{m}\right)
\end{equation}
with reference to a player's deck of size $m$ with $k$ special cards and separations $\ell_0,\ell_1,\dots,\ell_k$, in the notation of \eqref{eq:intervalsSEP}-\eqref{eq:ellsSEP} above.
For short, $H = -\sum_{j=0}^{k} \frac{\ell_j}{m} \ln\left(\frac{\ell_j}{m}\right)$ with the standard convention that $0 \ln(0) = 0$.

During the game the position entropy is constrained as
\begin{equation}
H_{\min} \,\leqslant\, H \,\leqslant\, H_{\max}
\end{equation}
with
\begin{equation}\label{eq:HminHmaxPosEntr}
H_{\min} \,=\, -\left(1 - \frac{k}{m}\right) \ln\left(1 - \frac{k}{m}\right)\,, \qquad H_{\max} \,=\, \ln(k+1)\,.
\end{equation}
Maximum entropy $H_{\max}$ occurs when special cards partition the deck into equal intervals (uniform dispersion). Minimum entropy $H_{\min}$ occurs when all special cards are clustered at one boundary of the deck.

\begin{figure}[t!] 
    \begin{adjustwidth}{-4cm}{-3.7cm}
        \centering
        \begin{minipage}[b]{0.55\textwidth} 
            \centering
            \includegraphics[width=\linewidth]{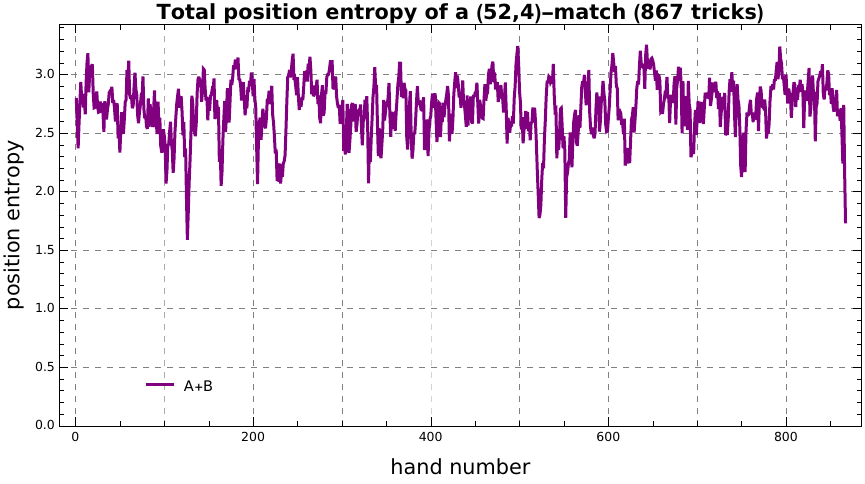}

            \includegraphics[width=\linewidth]{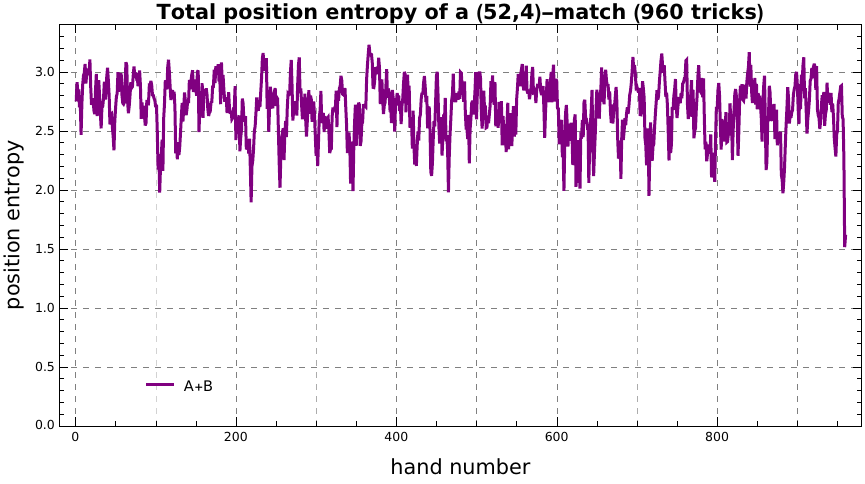}
        \vspace{0.1cm}
        \end{minipage}
        \quad 
        \begin{minipage}[b]{0.55\textwidth} 
            \centering
            \includegraphics[width=\linewidth]{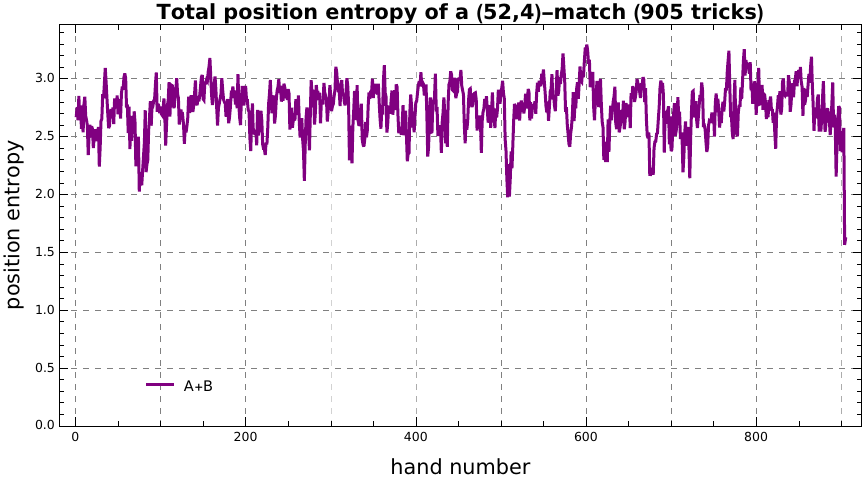}

            \includegraphics[width=\linewidth]{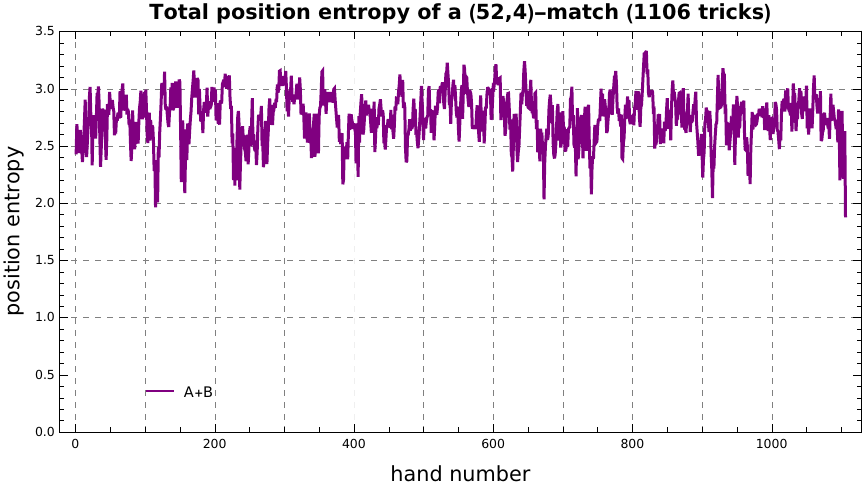}

          \hspace{1.2cm}
        \end{minipage}
        \caption{\small Total position entropies (sum of player $A$'s and player $B$'s position entropies) in the course of some ultra-long matches. \normalsize}
        \label{fig:position_entropy_figure}
    \end{adjustwidth}
\end{figure}

\begin{table}[!htbp]
\centering
\begin{tabular}{||c|c|c|c|c||}
\hline
$\quad N\quad$ & $\quad\mathfrak{R}\quad$ & \makecell{final decks' average  \\  position entropy $H_{\text{final}}$} & $\quad H_{\text{min}} \quad $ & $\quad H_{\text{max}} \quad$ \\
\hline
52 & 4 & 1.738 & 0.254 & 2.833 \\
\hline
52 & 3 & 1.706 & 0.202 & 2.565 \\
\hline
40 & 3 & 1.568 & 0.245 & 2.565 \\
\hline
40 & 2 & 1.373 & 0.179 & 2.197 \\
\hline
\end{tabular}
\vspace{0.2cm}
\caption{\small Position entropy in the winner's final deck averaged over randomly selected long $(N,\mathfrak{R})$-matches, and comparison with the theoretical upper and lower bounds \eqref{eq:HminHmaxPosEntr}. Here `long' corresponds to more than 350 tricks per match.\normalsize}\label{tab:positionentropy}
\end{table}

 The evolution of total position entropies in the course of ultra-long matches is displayed in Figure \ref{fig:position_entropy_figure}. For concreteness, in the standard $(52, 4)$-setting, where the winner ends up with a complete deck of $m = N = 52$ cards and $k = 4\mathfrak{R} = 16$ special cards, the bounds on position entropy are
\begin{equation*}
H_{\text{min}} \,\approx\, 0.255\,, \qquad H_{\text{max}} \,\approx\, 2.833\,.
\end{equation*}
 However, empirical analysis of ultra-long matches (Figure \ref{fig:separation}) reveals that the final total position entropy consistently lies in the range
\begin{equation*}
1.5 \,\lesssim\, H_{\text{final}} \,\lesssim\, 2.0\,,
\end{equation*}
indicating that the winner's deck exhibits moderate dispersion of special cards at match conclusion, with no accumulation of all special cards in tight clusters (which would yield $H_{\text{final}} \approx H_{\text{min}}$). Similar findings for other $(N,\mathfrak{R})$-settings are reported in Table \ref{tab:positionentropy}.

 Even in winning positions, the deck maintains significant structural inhomogeneity, with special cards neither maximally clustered nor uniformly distributed. Combined with the separation analysis (Figure \ref{fig:separation}), this provides evidence of \emph{persistent partial randomisation} of special card positions throughout the match: the game dynamics continuously redistribute special cards sufficiently to prevent extreme clustering, but do not drive the system toward perfect uniform dispersion. This intermediate regime -- sustaining moderate entropy values throughout hundreds of tricks -- characterises the quasi-equilibrium phase observed in ultra-long matches.

\section{Existence of infinite games}\label{sec:R1}

Non-terminating matches for Beggar-My-Neighbour do, in fact, exist. Since the total number of deck configurations is finite, these are games eventually consisting of a repeating loop of tricks.

At a sufficiently low ratio of special to total cards, $4\mathfrak{R}/N$, extensive numerical simulations are more likely to explicitly identify actual infinite matches, if any. In the setting $(N,\mathfrak{R})=(40,1)$ ($4\mathfrak{R}/N= 10\%$), for example, we successfully identified an infinite match starting with the initial decks
\begin{equation}\label{eq:inf01}
\begin{split}
\texttt{DeckA} &= \texttt{[CCCCCCCCCCC1CCCCCC1C]}, \\
\texttt{DeckB} &= \texttt{[1CCCCCCCCCCCCC1CCCCC]}
\end{split}
\end{equation}
 in a simulation of $10^7$ distinct initial configurations. For larger values of the ratio $4\mathfrak{R}/N$, including the standard $(40,3)$-scenario ($4\mathfrak{R}/N= 30\%$) and $(52,4)$-scenario ($4\mathfrak{R}/N\approx 31\%$), we have already argued that exploring sample sizes of this magnitude is probabilistically insufficient to detect infinite matches. For instance, as noted in Section \ref{sec:num}, no infinite $(52,4)$-matches were detected or reported by \cite{Paulhus_Beggar1999,Casella2024} in brute-force simulations ranging from $10^9$ to $10^{15}$ distinct decks. Only recently, using a general scheme for the reversed reconstruction of loops, the work \cite{Casella2024} exhibited a family of initial $(52,4)$-decks, including
\begin{equation}\label{eq:inf02}
\begin{split}
\texttt{DeckA} &= \texttt{[CCC3CCC2C3241CCCCC441CC1CC]}, \\
\texttt{DeckB} &= \texttt{[CCCCCCCCCC2CCCC32C1CCCCC34]},
\end{split}
\end{equation}
all of which evolve into the same loop, thereby yielding an infinite match.

We note that, presumably due to a chain of misreporting in the literature, it has frequently been claimed (for instance in \cite{Berlekamp-Conway-Guy-2004,Spivey2010,Lakshtanov-Aleksenko-Beggar2013,Beggar-Gentilini}) that infinite (52,4)-matches were found in \cite{Paulhus_Beggar1999} with a frequency of approximately 1 in 150,000. Both assertions are false. In reality, \cite{Paulhus_Beggar1999} identified loops only in smaller, modified deck configurations, such as $N=24$ with two -- rather than four -- special cards for each type 1, 2, 3, 4 (corresponding to J, Q, K, A, as noted in Section \ref{sec:game_rules}).

\begin{figure}[!htbp]
    \centering

    \begin{minipage}{1.0\textwidth}
        \centering
        \includegraphics[width=0.8\textwidth]{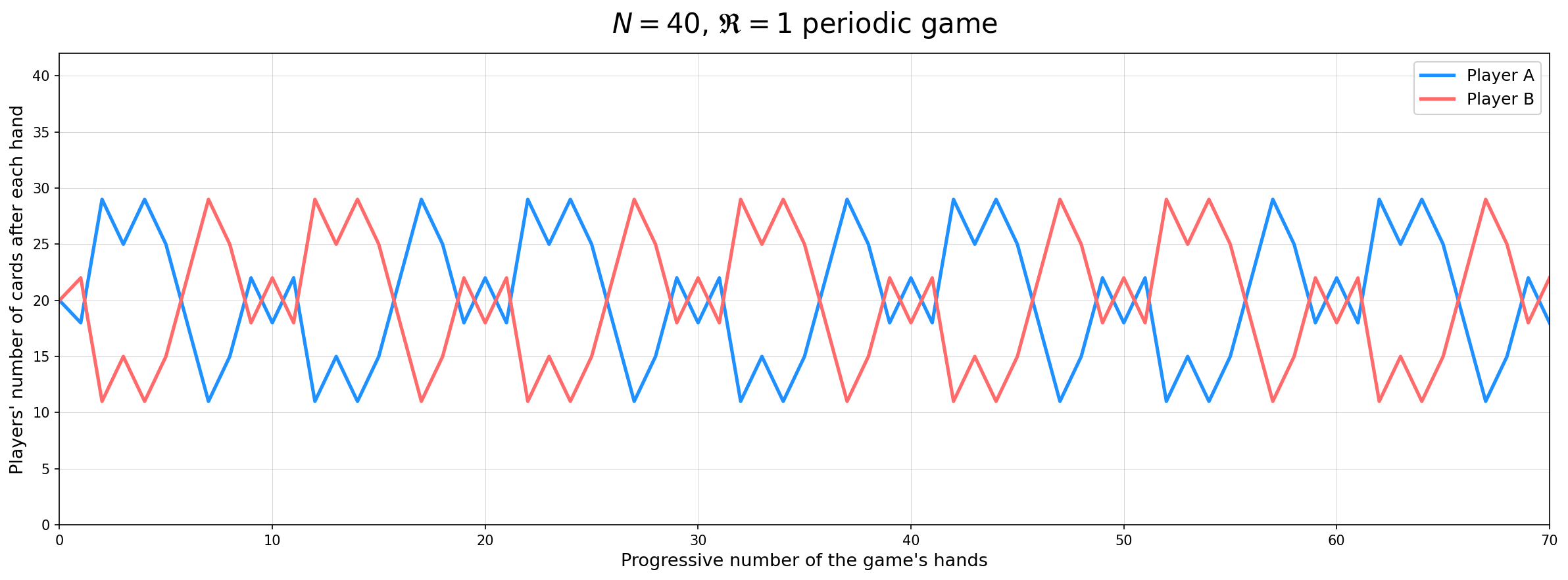}
    \end{minipage}

    \vspace{0.5cm} 

    \begin{minipage}{0.48\textwidth}
        \centering
        \includegraphics[width=\linewidth]{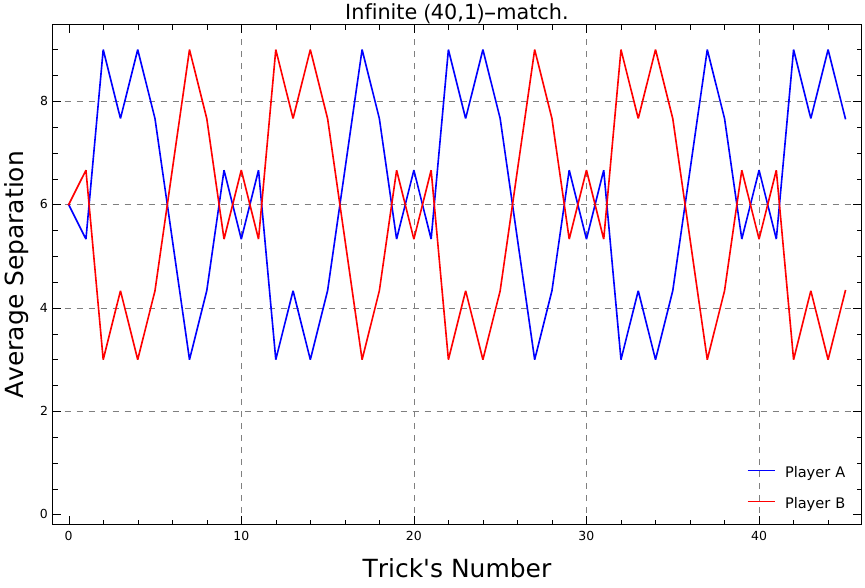}
    \end{minipage}
    \hfill 
    \begin{minipage}{0.48\textwidth}
        \centering
        \includegraphics[width=\linewidth]{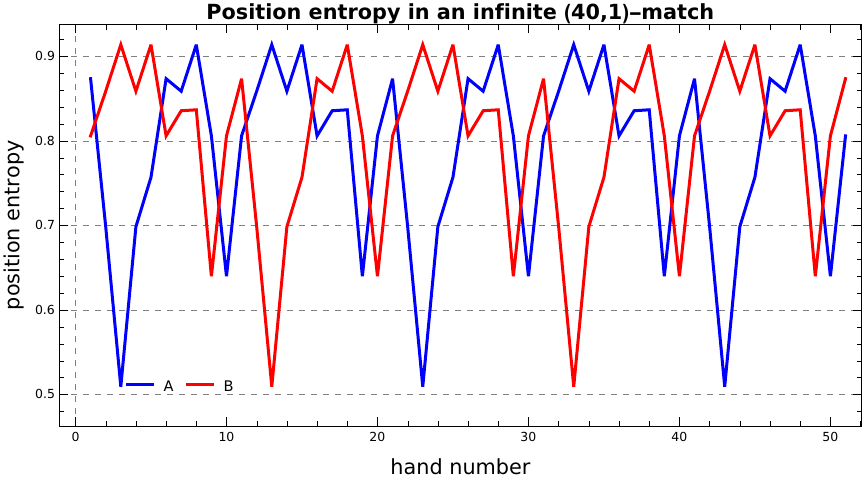}
    \end{minipage}

    \caption{\small Dynamics of the infinite $(40,1)$-match with initial decks \eqref{eq:inf01}: decks's number of cards (top), average separation of rank cards (bottom left), position entropies (bottom right). The loop begins at the first trick, with each player's deck oscillating between 11 and 29 cards (20-trick period). Periodicity in deck composition was verified computationally. \normalsize}
    \label{fig:40_1_match_combined}
\end{figure}
\begin{figure}[!htbp]
    \centering

    \begin{minipage}{1.0\textwidth}
        \centering
        \includegraphics[width=0.8\textwidth]{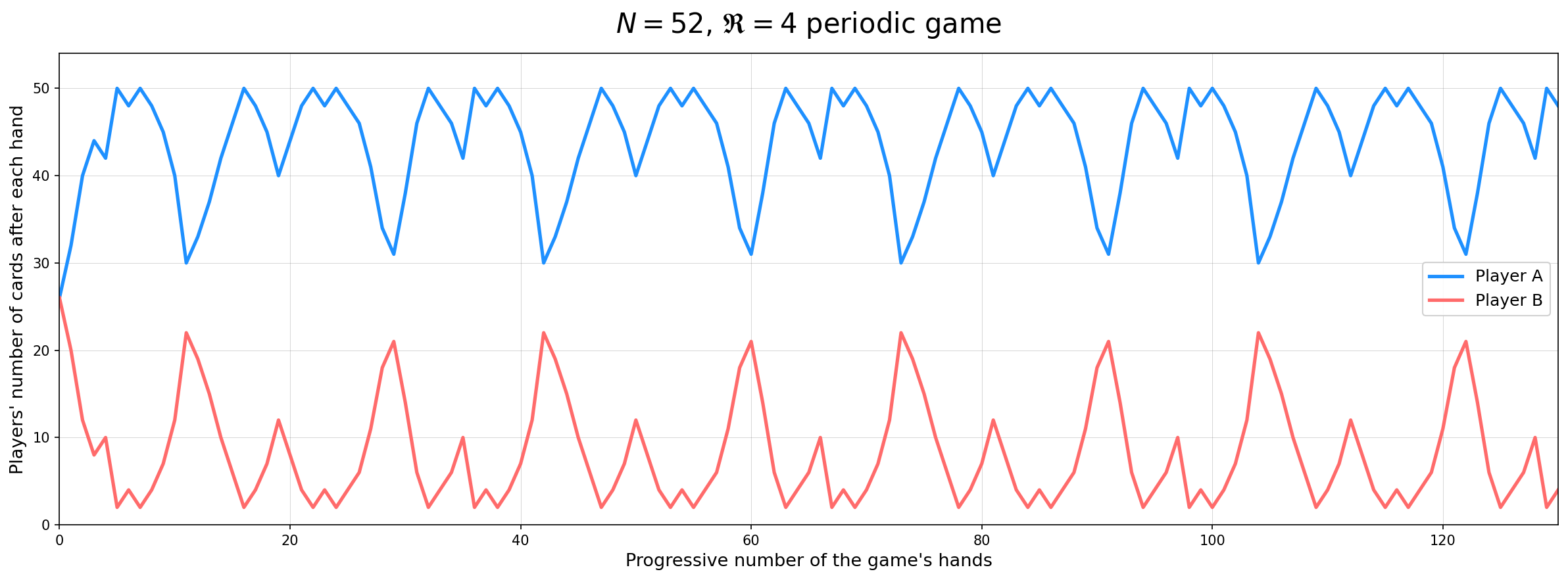}
    \end{minipage}

    \vspace{0.5cm} 

    \begin{minipage}{0.48\textwidth}
        \centering
        \includegraphics[width=\linewidth]{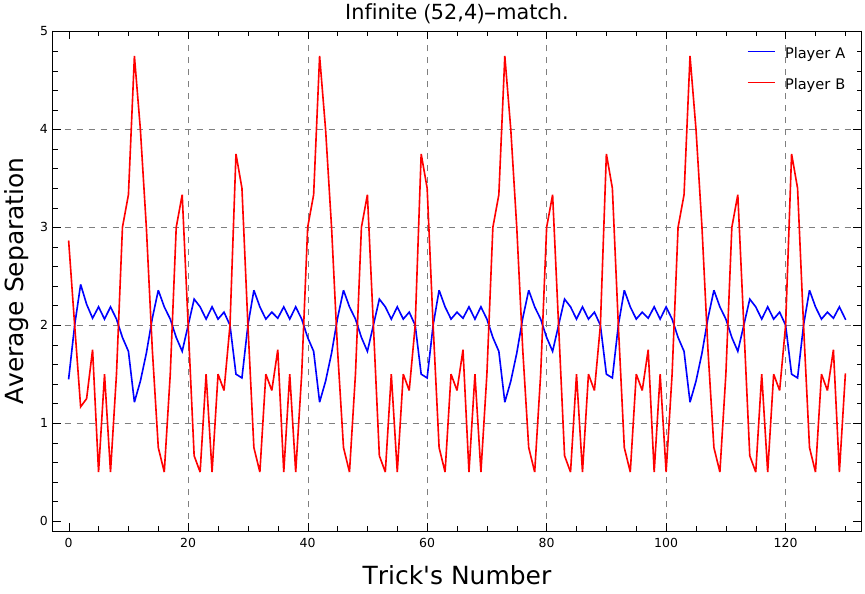}
    \end{minipage}
    \hfill 
    \begin{minipage}{0.48\textwidth}
        \centering
        \includegraphics[width=\linewidth]{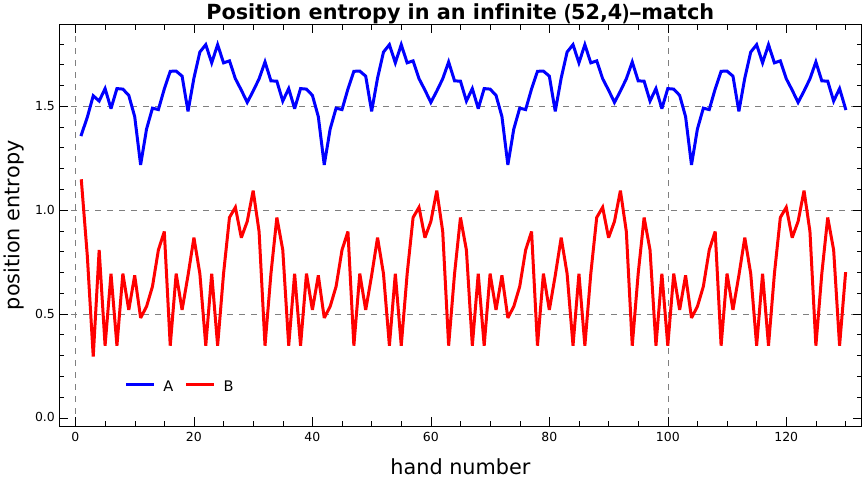}
    \end{minipage}

    \caption{\small Dynamics of the infinite $(52,4)$-match with initial decks \eqref{eq:inf02}: decks's number of cards (top), average separation of rank cards (bottom left), position entropies (bottom right). The loop begins after 4 tricks. Player $A$'s deck oscillates between 30 and 50 cards (62-trick period), maintaining a size strictly larger than player $B$'s deck, which oscillates between 2 and 22 cards. Periodicity in deck composition was verified computationally. \normalsize}
    \label{fig:52_4_match_combined}
\end{figure}

Relevant features of the infinite matches beginning with configurations \eqref{eq:inf01} and \eqref{eq:inf02} are shown in Figures \ref{fig:40_1_match_combined} and \ref{fig:52_4_match_combined}.

We observe that in either infinite match the system never returns to an exact $50\%$--$50\%$ distribution at any subsequent integer trick within the loop.

The considered $(40,1)$-loop displays an exact inversion of roles between the two players at the half-period mark, regarding card counts, average separation of special cards, and position entropy. The $20$-trick cycle consists of two consecutive $10$-trick phases, each characterised by a dominant player who holds more cards and exhibits lower average separation (more clustered special cards). At the transition, the player who was dominant (large deck, clustered specials) becomes subordinate (small deck, dispersed specials), and vice versa. The return to the initial balanced state closes the cycle.

It is plausible that this $(40,1)$-loop occupies a special, highly constrained region of the state space (of size $\sim 10^{14}$). Random sampling is unlikely to encounter configurations with exact half-period swap symmetry; the successful identification of this cycle via forward simulation ($10^7$ trials) may have been fortuitous, or it may indicate that cycles with high symmetry occupy disproportionately large basins of attraction relative to their measure-theoretic volume in state space. This hypothesis merits further investigation via systematic backward search from the cycle states.

The considered $(52,4)$-loop exhibits persistent asymmetry in the deck sizes, with player $A$ consistently maintaining a larger deck than player $B$, with a deck size difference ranging from 8 to 48 cards. The two players adopt fundamentally different roles within the cycle: player $A$ acts as a `reservoir', retaining the majority of cards and sustaining the capacity to initiate challenges, while player $B$ operates with a volatile, smaller deck that periodically contracts to near-depletion (as few as 2 cards) before recovering through trick wins. Even when player $B$'s deck repeatedly approaches zero, this does not cause termination, as player $B$ retains the special cards necessary to respond to challenges and trigger the card-transfer dynamics that restore balance within the cycle.

In such a match, the two players' position entropies are comparable to the entropy profiles observed in ultra-long terminating matches before eventual termination (Section \ref{sec:et}, Figure \ref{fig:position_entropy_figure}). In this respect, ultra-long terminating matches may represent trajectories that pass close to the basins of attraction of infinite cycles without entering them.

 \section{Constructing initial decks for non-terminating matches}\label{sec:constructingInf}

The recent work \cite{Casella2024} presented a constructive approach to the identification of non-terminating matches -- specifically, in the $(52,4)$-setting -- consisting of two main phases.
\begin{enumerate}
 \item[1.] Construction of two non-balanced decks (i.e., with cards not split evenly) yielding a loop. The players' decks are assembled
 \begin{itemize}
  \item[1.1] by suitably concatenating elementary building blocks that separately, in a reduced non-four-suited $(N,\mathfrak{R})$-version of the game, are checked by direct inspection to produce a recurring cycle;
  \item[1.2] and by suitably inserting additional special cards, up to fulfilling the decks' composition for the $(52,4)$-setting, through a trial-and-error process of individual insertions, deletions, or substitutions, and subsequent tests for termination.
 \end{itemize}
 \item[2.] Backward dynamics reconstruction from certain configurations of the loop, until a balanced $26$--$26$ deck is obtained.
\end{enumerate}
 Proceeding this way, \cite[Figure 3]{Casella2024} reports one loop, along which decks remain unbalanced, and 30 initial (balanced) decks (one of which is the configuration \eqref{eq:inf02}) whose dynamics enters the loop after 4, or 5, or 6 tricks.

 Various elements of the above strategy were already present across the above-mentioned spectrum of amateur literature: \cite{Casella2024} has the valuable merit to have organised and implemented them coherently and successfully, so as to produce the first explicit non-terminating game.

 On the other hand, for the purpose of reproducing, systematising, and generalising such an approach, the way some of the steps of the reasoning of \cite{Casella2024} are presented
 appears as a somewhat intuitive, heuristic first step, even after an accurate reading, and some other steps deserve to be supplemented with clarifying remarks.

\begin{figure}[!htbp]
    \begin{adjustwidth}{-2cm}{-2cm}
        \centering

        \includegraphics[height=1.5cm]{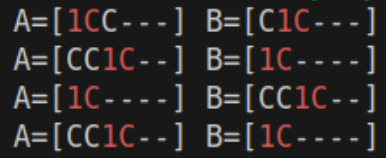}
        \;
        \includegraphics[height=1.5cm]{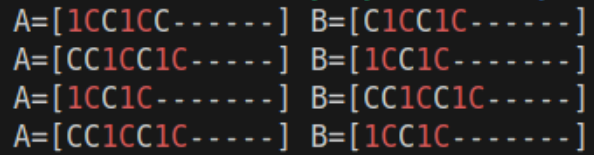}
        \;
        \includegraphics[height=1.5cm]{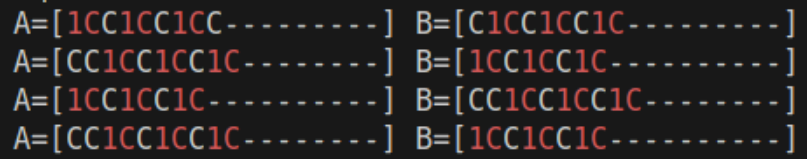}

        \caption{\small First three tricks from the initial decks \eqref{eq:inf03} (left), \eqref{eq:inf04} (centre), and \eqref{eq:inf05} (right). The third trick produces a configuration identical to the one emerging after the first trick and the dynamics enter a 2-tricks-long loop. \normalsize}
        \label{fig:three_basic_decks}
    \end{adjustwidth}
\end{figure}

 The intuition behind part 1.1 is surely the key and clever one: to start with minimal settings and initial configurations that by direct and easy inspection are seen to give rise to a periodic dynamics. \cite{Casella2024} starts with the following configuration, well familiar in the amateur literature,
\begin{equation}\label{eq:inf03}
\begin{split}
\texttt{DeckA} &= \texttt{[1CC]}, \\
\texttt{DeckB} &= \texttt{[C1C]},
\end{split}
\end{equation}
 and also with its doubling
  \begin{equation}\label{eq:inf04}
\begin{split}
\texttt{DeckA} &= \texttt{[1CC1CC]}, \\
\texttt{DeckB} &= \texttt{[C1CC1C]}.
\end{split}
\end{equation}
 The reasoning can be iterated, for example we may consider
   \begin{equation}\label{eq:inf05}
\begin{split}
\texttt{DeckA} &= \texttt{[1CC1CC1CC]}, \\
\texttt{DeckB} &= \texttt{[C1CC1CC1C]},
\end{split}
\end{equation}
etc. In either case, after the first trick, the match enters a 2-tricks-long loop (Figure \ref{fig:three_basic_decks}). In practice, the trick rule function $\mathfrak{F}$ `factorises' over such modular decks, as if they behave independently.

%

Then \cite{Casella2024} reports a subsequent procedure of insertions, deletions, or substitutions of additional ordinary and special cards, followed by tests for termination, carried out with a maximum of three simultaneous operations per trial. This process (part 1.2 of the summary outlined at the beginning of this Section) appears to be solely driven by trial-and-error modifications, exhaustive search, and filtering for non-termination and eventual standard composition of two decks in the four-suited $(52,4)$-setting. \cite{Casella2024} does not provide further computational details or summary statistics of the search process.

This eventually led to the identification of two configurations of unbalanced decks yielding, respectively, a 72-tricks loop and a 62-tricks loop.

%
%

 A crucial observation is that along such loops player~$A$'s deck can be viewed as decomposed into blocks terminated by rank-1 cards, which rotate through the deck over intervals of 4 tricks, while player~$B$'s deck periodically returns to the simple configuration \texttt{[1,C]}. Each rank-1-terminated block in player $A$’s deck contains exactly one rank-1 card at
 its end.

 For example, an initial configuration found for the 72-tricks loop is
 \begin{equation*}
\begin{split}
\texttt{DeckA} &= \texttt{[}\underbrace{\texttt{C\dots C1}}_{\text{20 C's, then 1}}\underbrace{\texttt{CC2C2C2C2C3C31}}_{\text{block 2}}\;\;\underbrace{\texttt{CC3C3C4C4C4C4C1}}_{\text{block 3}}\;\;\texttt{C]}\,, \\
\texttt{DeckB} &= \texttt{[1C]}
\end{split}
\end{equation*}
 Along the dynamics, between the checkpoints at 4, 8, 12, .... tricks, player~B's deck temporarily varies,
\begin{equation*}
\begin{split}
\text{trick } 0: \quad & \texttt{DeckB} = \texttt{[}\boxed{\texttt{1C}}\texttt{]} \\
\text{trick } 1: \quad & \texttt{DeckB} = \texttt{[CC}\boxed{\texttt{1C}}\texttt{]} \quad \text{(won 2 cards)} \\
\text{trick } 2: \quad & \texttt{DeckB} = \texttt{[CCCCC}\boxed{\texttt{1C}}\texttt{]} \quad \text{(won 4 more cards)} \\
\text{trick } 3: \quad & \texttt{DeckB} = \texttt{[CCCCCCCCCCC}\boxed{\texttt{1C}}\texttt{]} \quad \text{(won 6 more cards)} \\
\text{trick } 4: \quad & \texttt{DeckB} = \texttt{[}\boxed{\texttt{1C}}\texttt{]} \quad \text{(reset after challenge resolution)},
\end{split}
\end{equation*}
and meanwhile one of player $A$’s rank-1-terminated blocks cycles from near the top to near the bottom:
\begin{equation*}
\begin{split}
\text{trick}\;\;\; 0: \quad & \texttt{DeckA} = \texttt{[}\boxed{\texttt{C\dots C1}}\;\;\texttt{CC2C2C2C2C3C31}\;\;\texttt{CC3C3C4C4C4C4C1}\;\;\texttt{C]} \\
\text{trick}\;\;\; 4: \quad & \texttt{DeckA} = \texttt{[CC2C2C2C2C3C31}\;\;\texttt{CC3C3C4C4C4C4C1}\;\;\boxed{\texttt{C\dots C1}}\;\;\texttt{C]} \\
\text{trick}\;\;\; 8: \quad & \texttt{DeckA} = \texttt{[CC3C3C4C4C4C4C1}\;\;\boxed{\texttt{C\dots C1}}\;\;\texttt{CC2C2C3C2C3C21}\;\;\texttt{C]} \\
\text{trick} \;\;12: \quad & \texttt{DeckA} = \texttt{[}\boxed{\texttt{C\dots C1}}\;\;\texttt{CC2C2C3C2C3C21}\;\;\texttt{CC4C3C4C3C4C41}\;\;\texttt{C]}
\end{split}
\end{equation*}
etc.

Such rank-1-terminated structure with \texttt{[1,C]}-resets was used in \cite{Casella2024} as heuristics that guided which small non-terminating games to use as building blocks (those with simple loops), how to concatenate them (aligning rank-1 terminators), what configuration to assign player $B$ (minimal \texttt{[1,C]} that enables resets), and how to organise the subsequent trial-and-error modifications in step~1.2 (preserving the rank-1-terminated structure while inserting higher-rank special cards).

The reasoning of \cite{Casella2024} rested on the twofold informal assumption of block independence (each rank-1-terminated block processes independently, without interfering with other blocks) and reliable resets (player $B$'s deck reliably returns to a simple configuration, e.g., \texttt{[1,C]} between blocks).

 These assumptions are only approximately satisfied. Strictly speaking, blocks are not independent: during each 4-tricks interval, player $B$'s deck varies significantly (\texttt{[1C]} $\to$ \texttt{[CC\;1C]} $\to$ \texttt{[CCCCC\;1C]} $\to$ \texttt{[CCCCCCCCCCC\;1C]} $\to$ \texttt{[1C]}), showing complex card exchanges. Moreover, the special cards (ranks 2, 3, 4) within player $A$'s deck undergo internal permutations and rearrangements during processing, indicating coupling between blocks that the simplified description of \cite{Casella2024} does not account for.

%

\begin{figure}[!htbp]
    \centering

    \includegraphics[width=0.9\textwidth]{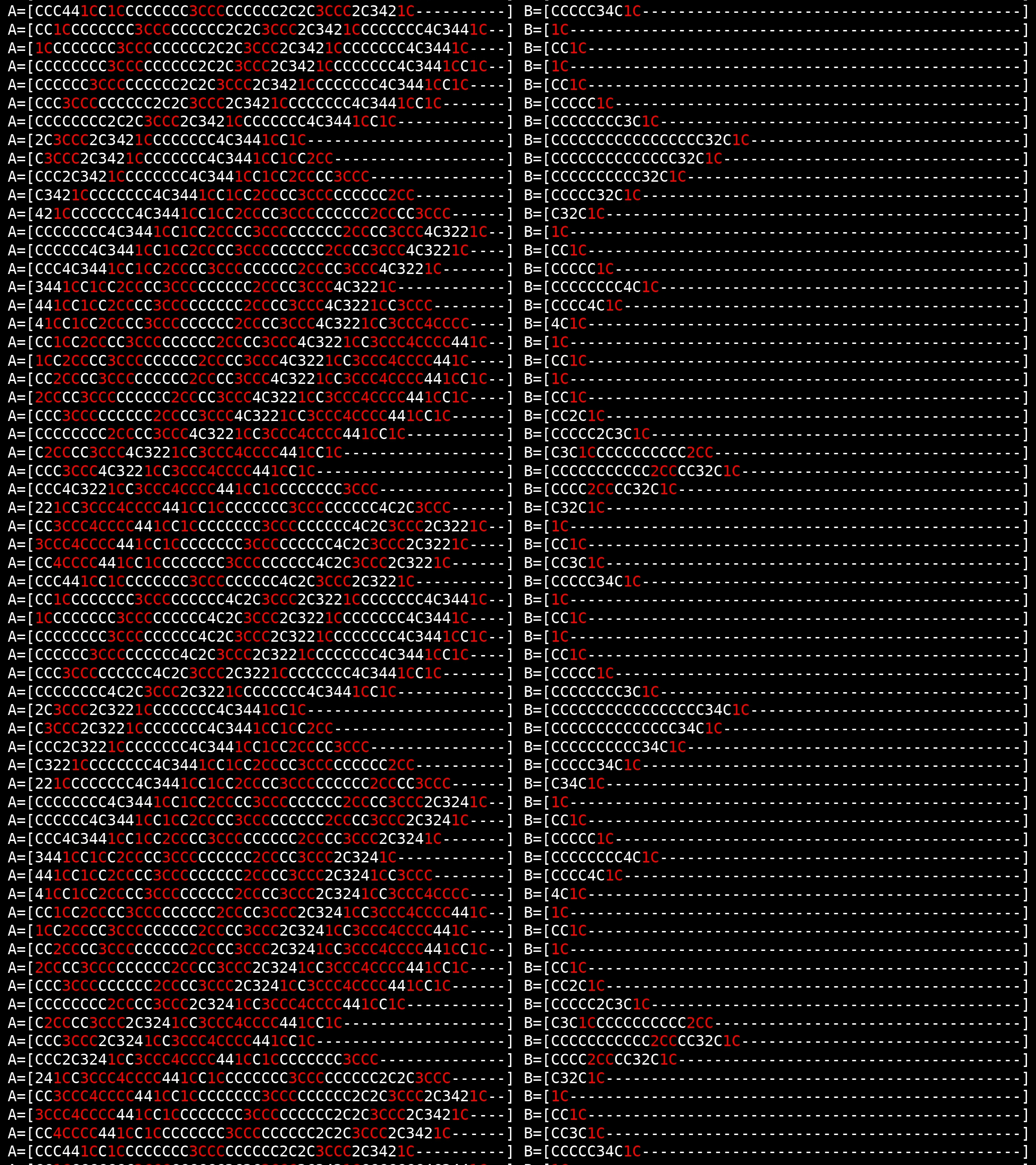}

    \caption{\small Dynamics of the 62-tricks loop identified in \cite{Casella2024}. \normalsize}
    \label{fig:casellaloop62}
\end{figure}

 The next phase (part 2) of the strategy outlined at the beginning of this Section concerns the quest for balanced decks whose dynamics reaches the loop. That is, the check that some of the configurations of the loop do belong to the reachable space \eqref{eq:reachablespace}.
 It appears that in \cite{Casella2024} this check was successful for the 62-tricks loop. This is precisely the loop analysed in Figure \ref{fig:52_4_match_combined}, whose dynamics is reported in Figure \ref{fig:casellaloop62}.


 As far as can be understood from the description in \cite{Casella2024}, this proceeded by
\begin{itemize}
 \item searching for new deck blocks with better ``backwards-play properties'',
 \item concatenating these blocks with rank-1 terminators,
 \item recognising that such concatenated decks enter a loop (by computational testing),
 \item applying exhaustive backward reconstruction (within computational limits) from states in this loop, building a tree of possible predecessor states, and searching for balanced decks among the predecessors.
\end{itemize}
``Backwards-play properties'' are informally described as the occurrence of independent behaviour of concatenated blocks and the preservation of non-termination when swapping special cards (which helps in achieving the four-suited $(52,4)$ composition). Recall also that in general backward play is non-deterministic (Section~\ref{sec:fwd-bkw-determinism}): given a state, there may be multiple predecessor states that lead to it.

Overall, the work \cite{Casella2024} is a proof of existence by exhibition of an explicit example, combined with methodological insights (blocks, templates, backward play), but not a complete algorithm. In particular, it does not provide a reproducible procedure connecting the 72-tricks loop (or any starting point) to the final balanced non-terminating configuration looping over 62 tricks. It describes the ideas that guided the search, but does not provide the step-by-step procedure that would allow one to derive the same solution deterministically.

\section{Infinite loops factory}\label{sec:loop_factory}

Inspired by the general, yet not fully algorithmic methodology of \cite{Casella2024} -- namely, the reconstruction of possibly unbalanced, four-suited $(N,\mathfrak{R})$-decks yielding non-terminating loops, carried out from small cycling decks via successive insertion trials -- we developed a C-script implementing systematic card insertion strategies with adaptive backtracking.

This aims to produce loops in various $(N,\mathfrak{R})$-settings, from whose states we exhaustively develop a backward reconstruction of pre-loop histories, until possibly hitting a balanced pair of decks.

At each step, the script tests whether inserting specific card
sequences preserves the non-termination property; successful insertions are saved as
checkpoints for later exploration. When a branch fails to progress, the algorithm backtracks
to a previously saved checkpoint -- selected either by proximity to target composition (via a
scoring function) or through pseudo-random selection to escape local minima -- and explores
alternative insertion strategies.

The search maintains up to 1500 checkpoints simultaneously, each recording a partial deck configuration along with its exploration depth, compositional score, and the strategy that generated it. This checkpoint buffer enables the algorithm to maintain multiple exploration paths in parallel, periodically pruning the oldest or least promising states when capacity is reached. The script may restart from the initial seed configuration to prevent permanent entrapment in unproductive regions of the search space.

The script implements distinct high-level strategies for determining which card types to prioritise during insertion.
\begin{itemize}
    \item \textbf{Low-rank-first}: prioritises filling special cards in ascending rank order (rank-1, then rank-2, \dots, up to rank-$\mathfrak{R}$), followed by ordinary cards. This strategy reflects the observation of \cite{Casella2024} that lower ranks often admit more positional flexibility.

    \item \textbf{High-rank-first}: inverse approach, attempting to insert rank-$\mathfrak{R}$ cards first, descending to rank-1, and concluding with ordinary cards. Useful when higher-rank cards require specific structural positioning.

    \item \textbf{Ordinaries-first}: fills all ordinary cards before attempting any special card insertions. This establishes a `skeleton' of common cards that may better accommodate subsequent special card placement.

    \item \textbf{Specials-first}: inserts all special cards (ranks 1 through $\mathfrak{R}$) before adding ordinary cards, creating a special-card framework that ordinary cards must then fit around.

    \item \textbf{Balanced fill}: attempts to maintain proportional progress toward the target composition, adding cards of each type in ratio roughly equal to their representation in the final $(N,\mathfrak{R})$-deck.

    \item \textbf{Cluster specials}: inserts special cards in concentrated groups, creating local `special-dense' regions within the deck.

    \item \textbf{Disperse specials}: spreads special cards with substantial ordinary-card buffers between them, minimising special-to-special adjacencies.

    \item \textbf{Rank-1 blocks}: constructs rank-1-terminated blocks following the template methodology of \cite{Casella2024}, building sequences of the form \texttt{[CC\ldots CC1]} as fundamental structural units.
\end{itemize}

For each high-level strategy, the script generates concrete card sequences using distinct patterns, varying from single-card insertions to burst insertions of up to ten cards:
\begin{itemize}
    \item \textbf{single}, individual cards \texttt{[C]}, \texttt{[1]}, \texttt{[2]}, \dots, \texttt{[$\mathfrak{R}$]}, prioritising special cards over ordinaries based on current deficiencies;
    \item \textbf{pairs}, two-card sequences including homogeneous pairs \texttt{[CC]}, \texttt{[$r$$r$]} and heterogeneous combinations \texttt{[C$r$]}, \texttt{[$r$C]} for all ranks $r \in \{1,\ldots,\mathfrak{R}\}$;
    \item \textbf{triplets}, three-card sequences with various special card placements \texttt{[CCC]}, \texttt{[C$r$C]}, \texttt{[CC$r$]}, \texttt{[$r$CC]};
    \item \textbf{short bursts (5 cards)}, sequences \texttt{[CCCCC]} for rapid ordinary-card addition;
    \item \textbf{medium bursts (7 cards)}, sequences \texttt{[CCCCCCC]} providing moderate-scale ordinary insertion;
    \item \textbf{long bursts (10 cards)},  maximal ordinary sequences \texttt{[CCCCCCCCCC]} for aggressive deck expansion;
    \item \textbf{interleaved}, alternating patterns \texttt{[C$r$C$r$C]} that distribute special cards within ordinary buffers;
    \item \textbf{rank-terminated}, sequences terminating in a special card, following the block construction of \cite{Casella2024}, namely, \texttt{[CC$r$]}, \texttt{[CCC$r$]}, \texttt{[CCCC$r$]};
    \item \textbf{mixed short bursts (up to 7 cards)}, sequences embedding a special card within ordinary buffers: \texttt{[CCC$r$CCC]}, \texttt{[C$r$C]}, \texttt{[CC$r$CC]}; this also supports multiple special cards \texttt{[C$r_1$$r_1$C]} when sufficient cards of rank $r_1$ remain needed;
    \item \textbf{mixed long bursts (up to 10 cards)}, extended mixed sequences \texttt{[CCCC$r$CCCCC]} and multi-special patterns \texttt{[CC$r_1$CC$r_2$CC]} incorporating two different ranks.
    \end{itemize}

Furthermore, the script employs bias modes controlling whether sequences are inserted into \texttt{deckA}, \texttt{deckB}, or both,
and implements corrective operations (card removal, card swapping) to handle compositional imbalances.
%
%

 A search depth is employed to measure the number of successful insertions from the initial seed configuration. Each successful insertion that preserves non-termination increments the depth counter and creates a checkpoint.

%
%

The checkpoint selection mechanism implements intelligent backtracking.
\begin{itemize}
 \item Checkpoints are scored by $N - \sum_{j=0}^{\mathfrak{R}} |\texttt{missing(}j\texttt{)}|$, where $\texttt{missing(}j\texttt{)}$ is a shorthand counting the deficiency or excess of cards of type $j$ ($j=0$ indicates ordinary cards). Higher scores indicate closer proximity to target composition.
 \item If the search returns to the same high-scoring checkpoint more than five consecutive times without progress, it pseudo-randomly selects an alternative checkpoint to escape the local minimum.
 \item Every 50 exploration rounds, the algorithm may restart from the initial seed to explore entirely different paths through the search space.
\end{itemize}

%
%
%
%

 The algorithm terminates upon finding a complete, non-terminating, $(N,\mathfrak{R})$-configuration, or alternatively, when the maximum attempt budget is exhausted, it reports the largest non-terminating configuration discovered, thus providing partial progress toward the target $(N,\mathfrak{R})$-deck.

 The loop factory produced a number of non-terminating configurations, with player $A$ always starting first, as summarised in Table \ref{tab:nonterminating_configs}.

\begin{table}[h]
\centering
\footnotesize 
\setlength{\tabcolsep}{4pt} 
\begin{tabular}{c@{\hskip 0.5cm}ll}
\hline
$(N,\mathfrak{R})$ & \texttt{DeckA} & \texttt{DeckB} \\
\hline
\multirow{2}{*}{$(12,1)$}
 & \texttt{[1CC1CC1CC]} & \texttt{[C1C]} \\
 & \texttt{[C1C1CC1CC]} & \texttt{[C1C]} \\
\hline
\multirow{3}{*}{$(20,2)$}
 & \texttt{[1CC]} & \texttt{[C2C2C1CC1CC2C2C1C]} \\
 & \texttt{[1CC2C2C1CC2C2C1CC]} & \texttt{[C1C]} \\
 & \texttt{[1CC1CC1CC2C]} & \texttt{[C2C2C1CC2]} \\
\hline
$(24,1)$ & \texttt{[1CCCCCCCC1CC1CCCCCCCC]} & \texttt{[C1C]} \\
\hline
$(24,2)$ & \texttt{[1CCCCC]} & \texttt{[CCCC1CC2C21CC2C21C]} \\
\hline
$(24,3)$ & \texttt{[1CC3C21CC2C33C32C21CC]} & \texttt{[C1C]} \\
\hline
$(28,1)$ & \texttt{[1CCCCCC1CCCCCC]} & \texttt{[CC1CCCCCCCCC1C]} \\
\hline
\multirow{3}{*}{$(28,3)$}
 & \texttt{[1CC3C2C1CC2C2C3C3C1CC2C]} & \texttt{[C3C1C]} \\
 & \texttt{[1CC1CC]} & \texttt{[C2C1CC3C3C2C3C2C3C2C1C]} \\
 & \texttt{[3C1CC3C1CC2C3C2C]} & \texttt{[C3C1CC2C2C1C]} \\
\hline
$(32,1)$ & \texttt{[1CCCCCCCC1CCCCCCCCCCCC1CCCCCC]} & \texttt{[C1C]} \\
\hline
$(32,2)$ & \texttt{[1CC2C1CC1CCCCCCCC2CCCC2CCCC]} & \texttt{[C1CC2]} \\
\hline
$(32,3)$ & \texttt{[2C1CC2C]} & \texttt{[C3C3C331CC2CCCC1CC2CCCC1C]} \\
\hline
\multirow{2}{*}{$(36,1)$}
 & \texttt{[1CC]} & \texttt{[CCCCCCC1CC1CCCCC1CCCCCCCCCCCCCCCC]} \\
 & \texttt{[1CCCCCCCCCCCCCC1CC1CCCCCCCCCCCCCC]} & \texttt{[C1C]} \\
\hline
\multirow{2}{*}{$(40,1)$}
 & \texttt{[1CCCCCCCCCCCCCCCCC1CCCCC1CCCCCCCC]} & \texttt{[CCCCC1C]} \\
 & \texttt{[1CCCCCCCC1CCCCCCCCCCCCCCCCC1CCCCCCCCC]} & \texttt{[C1C]} \\
\hline
$(52,1)$ & \texttt{[1CCCCC1CCCCCCCCCCCCCCCCC1CCCCCCCCCCCCCCCCCCCCCCCC]} & \texttt{[C1C]} \\
\hline
\end{tabular}
\medskip

\caption{\small Non-terminating $(N,\mathfrak{R})$-configurations discovered by the automated explorative search algorithm. Player $A$ always moves first.}
\label{tab:nonterminating_configs}
\end{table}

  For future searches of loop games, it is \emph{plausible to argue that loops exist only in an intermediate regime of the ratio $\mathfrak{R}/N$} -- that is (since decks are four-suited), an intermediate regime of $4\mathfrak{R}/N\in[0,1]$. This is suggested by the observation that both initial decks with no special cards ($4\mathfrak{R}/N=0\%$) and with all special cards ($4\mathfrak{R}/N=100\%$) yield finite games consisting of one trick only.

  Looking in retrospect at Table \ref{tab:nonterminating_configs}, as well as at \eqref{eq:inf02}, the evidence so far indicates that the regime of existence of loops at least includes values $4\mathfrak{R}/N$ between $7\%$ (scenario $(N,\mathfrak{R})=(52,1)$) and $43\%$ (scenario $(N,\mathfrak{R})=(28,3)$).

 \section{Infinite matches starting in a loop}\label{sec:inf-starting-loops}

\begin{figure}[!htbp]
    \centering

    \includegraphics[width=0.7\textwidth]{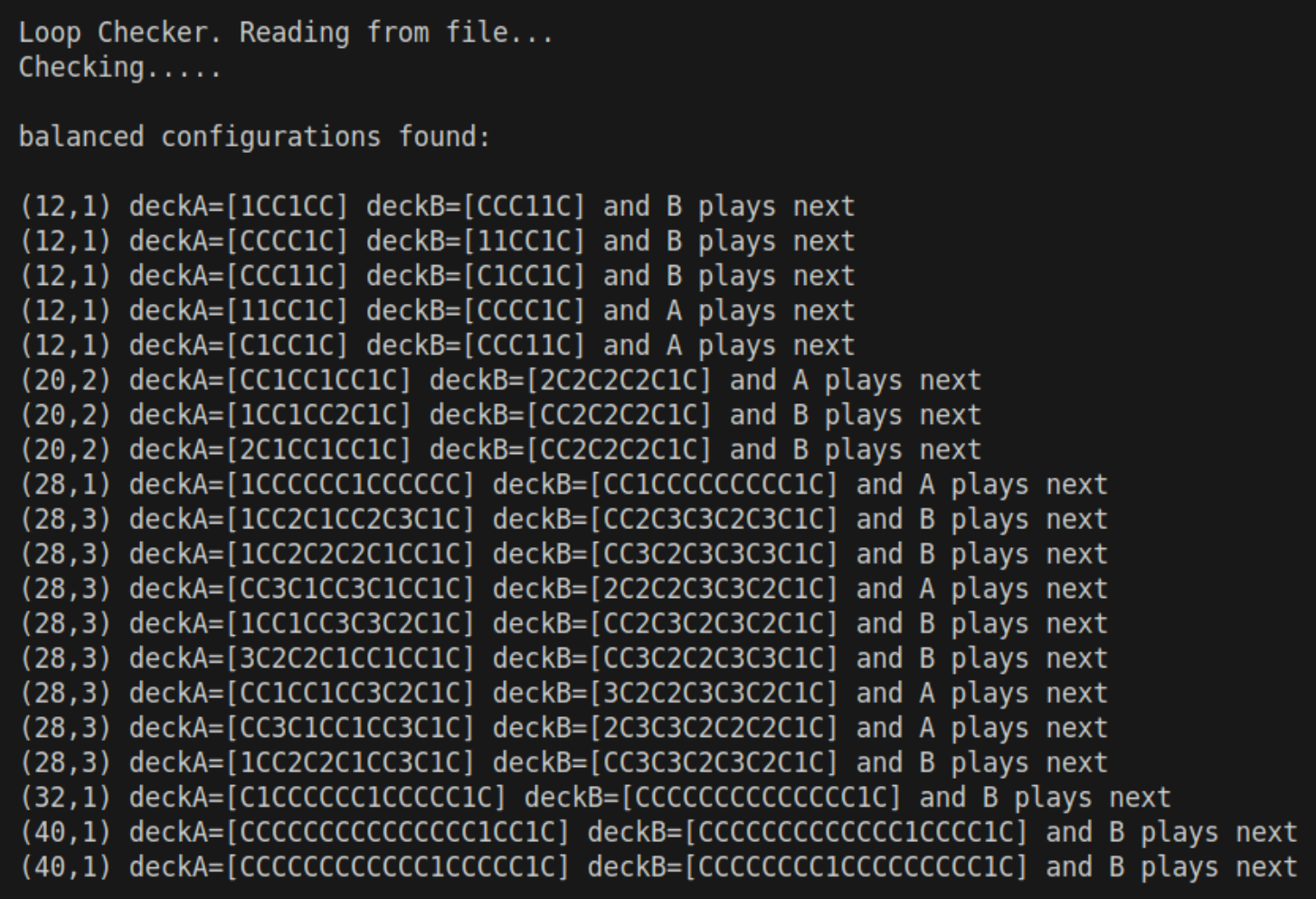}

    \caption{\small Identification of balanced configurations along the loops generated by the unbalanced decks of Table \ref{tab:nonterminating_configs}. \normalsize}
    \label{fig:loop_balance_check}
\end{figure}

 Investigation of the unbalanced, loop-generating configurations generated by our `loop factory' analysis (Section \ref{sec:loop_factory}) reveals that the vast majority of these loops transit through balanced decks (Figure \ref{fig:loop_balance_check}).

 This immediately yields initial decks for non-terminating games in various $(N,\mathfrak{R})$-settings. (In cases where Player $B$ has the lead in the identified balanced configuration, as occurs in several instances in Figure \ref{fig:loop_balance_check}, the roles of $A$ and $B$ are swapped for standardisation). They are reported in Table \ref{tab:balanced_loop_findings}.

\begin{table}[h]
\centering
\begin{tabular}{ccc@{\hskip 1cm}ll}
\hline
$(N,\mathfrak{R})$ & & $4\mathfrak{R}/N$ & \texttt{DeckA} & \texttt{DeckB} \\
\hline
\multirow{3}{*}{$(12,1)$} & & \multirow{3}{*}{$33.33\%$}
 & \texttt{[CCC11C]} & \texttt{[1CC1CC]} \\
 & & & \texttt{[11CC1C]} & \texttt{[CCCC1C]} \\
 & & & \texttt{[C1CC1C]} & \texttt{[CCC11C]} \\
\hline
\multirow{3}{*}{$(20,2)$} & & \multirow{3}{*}{$40.00\%$}
 & \texttt{[CC1CC1CC1C]} & \texttt{[2C2C2C2C1C]} \\
 & & & \texttt{[CC2C2C2C1C]} & \texttt{[1CC1CC2C1C]} \\
 & & & \texttt{[CC2C2C2C1C]} & \texttt{[2C1CC1CC1C]} \\
\hline
$(28,1)$ & & $14.29\%$ & \texttt{[1CCCCCC1CCCCCC]} & \texttt{[CC1CCCCCCCCC1C]} \\
\hline
\multirow{8}{*}{$(28,3)$} & & \multirow{8}{*}{$42.86\%$}
 & \texttt{[CC2C3C3C2C3C1C]} & \texttt{[1CC2C1CC2C3C1C]} \\
 & & & \texttt{[CC3C2C3C3C3C1C]} & \texttt{[1CC2C2C2C1CC1C]} \\
 & & & \texttt{[CC3C1CC3C1CC1C]} & \texttt{[2C2C2C3C3C2C1C]} \\
 & & & \texttt{[CC2C3C2C3C2C1C]} & \texttt{[1CC1CC3C3C2C1C]} \\
 & & & \texttt{[CC3C2C2C3C3C1C]} & \texttt{[3C2C2C1CC1CC1C]} \\
 & & & \texttt{[CC1CC1CC3C2C1C]} & \texttt{[3C2C2C3C3C2C1C]} \\
 & & & \texttt{[CC3C1CC1CC3C1C]} & \texttt{[2C3C3C2C2C2C1C]} \\
 & & & \texttt{[CC3C3C2C3C2C1C]} & \texttt{[1CC2C2C1CC3C1C]} \\
\hline
$(32,1)$ & & $12.50\%$ & \texttt{[CCCCCCCCCCCCCC1C]} & \texttt{[C1CCCCCC1CCCCC1C]} \\
\hline
\multirow{2}{*}{$(40,1)$} & & \multirow{2}{*}{$10.00\%$}
 & \texttt{[CCCCCCCCCCCCC1CCCC1C]} & \texttt{[CCCCCCCCCCCCCCC1CC1C]} \\
 & & & \texttt{[CCCCCCCC1CCCCCCCCC1C]} & \texttt{[CCCCCCCCCCCC1CCCCC1C]} \\
\hline
\end{tabular}
\medskip

\caption{\small Balanced initial decks ($|\texttt{DeckA}|=|\texttt{DeckB}|=N/2$) generating loops directly from the first trick. All configurations are standardised such that Player $A$ moves first.}
\label{tab:balanced_loop_findings}
\end{table}

The observation that a significant portion of computationally generated loops traverses balanced configurations -- in sharp contrast to the strictly unbalanced loop found in \cite{Casella2024} and discussed in Section \ref{sec:constructingInf}  -- can be attributed to the methodological differences in their discovery.

The constructive approach typically relies on rigid modular templates (e.g., rank-1 terminated blocks) that often lock players into fixed dynamical roles, such as a dominant `reservoir' versus a recurrently depleting `reset' deck, thereby structurally preserving asymmetry.

Conversely, the automated search algorithm described in Section \ref{sec:loop_factory} operates without such architectural constraints. By leveraging diverse insertion heuristics (particularly \textit{Balanced Fill} and \textit{Simultaneous} modes) combined with deep adaptive backtracking, the algorithm effectively samples the state space for cycles exhibiting high card-exchange mobility.

These naturally discovered loops allow the card count to oscillate freely across the $N/2$ equilibrium line, revealing the balanced entry points listed in Table \ref{tab:balanced_loop_findings} that purely constructive methods might overlook.

 \section{Dynamical landscape and open problems}\label{sec:retrosp-openpr}

 The multiple aspects analysed here of the Beggar-My-Neighbour card game reveal a rich dynamical landscape consisting of finite games whose length ranges from single tricks to ultra-long matches, as well as of infinite games entering a loop right at the beginning of the match or at some later trick.

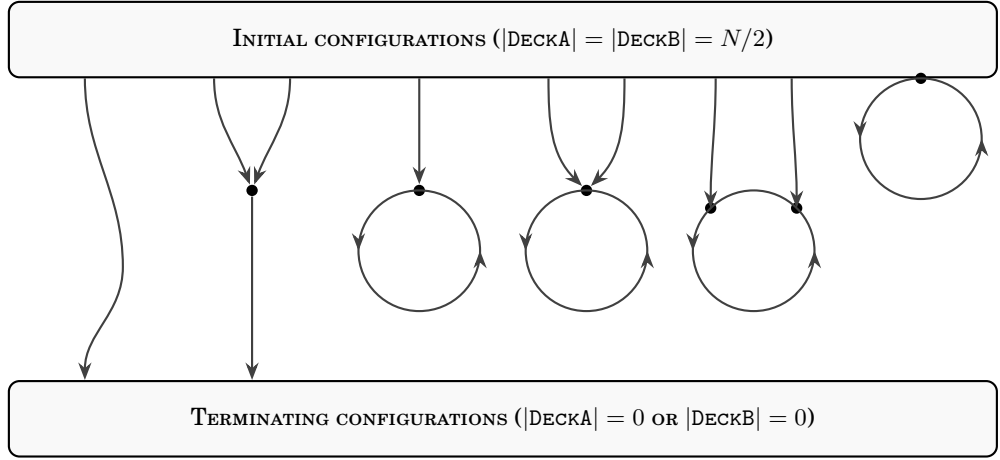
\begin{figure}[!htbp]
    \centering
    \begin{tikzpicture}[
        >=Stealth, 
        node distance=2cm,
        state_set/.style={
            draw,
            rectangle,
            rounded corners,
            minimum width=13cm, 
            minimum height=1cm,
            align=center,
            fill=gray!5,
            thick,
            font=\small\scshape
        },
        dot/.style={
            circle,
            fill=black,
            inner sep=1.5pt
        },
        flow/.style={
            ->,
            thick,
            darkgray,
            smooth,
            tension=0.7
        },
        doublearrow/.style={
            postaction={decorate},
            decoration={
                markings,
                mark=at position 0.25 with {\arrow{Stealth}},
                mark=at position 0.75 with {\arrow{Stealth}}
            }
        }
    ]

    \node[state_set] (start) {Initial configurations ($|\texttt{DeckA}|=|\texttt{DeckB}|=N/2$)};
    \node[state_set, below=4cm of start] (end) {Terminating configurations ($|\texttt{DeckA}|=0$ or $|\texttt{DeckB}|=0$)};


    \coordinate (a_top) at ($(start.south) + (-5.5,0)$);
    \coordinate (a_bot) at ($(end.north) + (-5.5,0)$);
    \draw[flow] (a_top) to[out=-90, in=90] ++(0.5,-2.5) to[out=-90, in=90] (a_bot);

    \coordinate (b1_top) at ($(start.south) + (-3.8,0)$);
    \coordinate (b2_top) at ($(start.south) + (-2.8,0)$);
    \coordinate (b_merge) at (-3.3, -2.0);
    \coordinate (b_bot) at ($(end.north) + (-3.3,0)$);
    \node[dot] (merge_point) at (b_merge) {};

    \draw[flow] (b1_top) to[out=-90, in=110] (merge_point);
    \draw[flow] (b2_top) to[out=-90, in=70] (merge_point);
    \draw[flow] (merge_point) to[out=-90, in=90] (b_bot);

    \coordinate (loop1_c) at (-1.1, -2.8);
    \node[dot] (loop1_entry) at ($(loop1_c)+(0, 0.8)$) {};
    \coordinate (c_top) at ($(start.south) + (-1.1,0)$);

    \draw[flow] (c_top) to[out=-90, in=90] (loop1_entry);
    \draw[thick, darkgray, doublearrow] (loop1_entry) arc (90:450:0.8);

    \coordinate (loop2_c) at (1.1, -2.8);
    \node[dot] (loop2_entry) at ($(loop2_c)+(0, 0.8)$) {};

    \coordinate (d1_top) at ($(start.south) + (0.6,0)$);
    \coordinate (d2_top) at ($(start.south) + (1.6,0)$);

    \draw[flow] (d1_top) to[out=-90, in=135] (loop2_entry);
    \draw[flow] (d2_top) to[out=-90, in=45] (loop2_entry);
    \draw[thick, darkgray, doublearrow] (loop2_entry) arc (90:450:0.8);

    \coordinate (loop3_c) at (3.3, -2.8);

    \coordinate (entry_e1) at ($(loop3_c) + (135:0.8)$);
    \node[dot] at (entry_e1) {};
    \coordinate (e1_top) at ($(start.south) + (2.8,0)$);
    \draw[flow] (e1_top) to[out=-90, in=90] (entry_e1);

    \coordinate (entry_e2) at ($(loop3_c) + (45:0.8)$);
    \node[dot] at (entry_e2) {};
    \coordinate (e2_top) at ($(start.south) + (3.8,0)$);
    \draw[flow] (e2_top) to[out=-90, in=90] (entry_e2);

    \draw[thick, darkgray, doublearrow] ($(loop3_c)+(0,0.8)$) arc (90:450:0.8);

    \coordinate (f_top) at ($(start.south) + (5.5,0)$);
    \node[dot] (f_entry) at (f_top) {};

    \draw[thick, darkgray, doublearrow] (f_entry) arc (90:450:0.8);

    \end{tikzpicture}
    \caption{\small Schematic representation of the game's dynamical landscape. From left to right: (1) unique terminating trajectory; (2) merging terminating trajectories; (3) single transient trajectory entering a loop; (4) multiple transient trajectories entering a loop at the same point; (5) separate transient trajectories entering a loop at distinct points; and (6) immediate loops where the balanced initial configuration is part of the cycle.}
    \label{fig:game_topology}
\end{figure}

  Figure \ref{fig:game_topology} depicts our understanding of this dynamical landscape, where we also reasonably conjecture the existence of merging terminating trajectories as well as separate transient trajectories entering the same loop.

 Brute-force numerical exploration of a tiny yet significant portion of the huge initial configurations space (Sections \ref{sec:num} and \ref{sec:et}) indicates an \emph{approximate} geometric regime for the distribution of terminating games in terms of their duration (more precisely, for the match duration $\mathcal{T}$ random variable), even though the dynamics of each individual game is perfectly deterministic, with no degree of stochasticity at all.

 The geometric (exponential) distribution approximation (Figure \ref{fig:semilog_40_3} and Section \ref{sec:et}) on the one hand provides an order-of-magnitude assessment of the probability of observing ultra-long matches (Section \ref{sec:et}), and on the other reasonably accounts for the somewhat erratic, stock-market-like dynamics of the decks’ cardinalities in the course of ultra-long matches (Figures \ref{fig:tricks_oscillation_plot_R3} and \ref{fig:tricks_oscillation_plot_52-4}). Approximately, the match duration $\mathcal{T}$ displays the typical memory-less property of geometric random variables (equations \eqref{eq:memory-less}-\eqref{eq:memory-less3}) and therefore, approximately, the fact that a game has endured for a certain amount of time has no predictive power over its future length -- in this respect we clarified the point of view of \cite{Casella2024}: each game is fully deterministic, hence the memory-less property should be understood as an approximation valid for the typical duration range, not as a property across all game lengths.

 In fact, due again to approximate geometric features, similar erratic oscillatory patterns were found for the evolution of the average separation among special cards and of position entropies (Section \ref{sec:et}, Figures \ref{fig:separation} and \ref{fig:position_entropy_figure}), with persistent partial randomisation of special card positions throughout the match preventing extreme clustering as well as perfect uniform dispersion.

 This explains why we couldn't really identify indicators of future duration consisting of simple measurements such as deck's cardinality, special cards' separation, or position entropies, and ultimately we could not identify meaningful metrics in the state space of the game -- a structure made clear in the abstract formalism laid down in Sections \ref{sec:abstr}-\ref{sec:fwd-bkw-determinism}.

 Finally, brute-force exploration combined with a refined reversed reconstruction of game trajectories reveal a rich scenario also for infinite games entering a loop at the beginning of the match or after some transient (Sections \ref{sec:R1}-\ref{sec:inf-starting-loops}).

 The above considerations bring us to a core of non-trivial, instructive open questions that certainly deserve to be addressed and are the object of our future investigation. \emph{Apparently} they are specific to the Beggar-My-Neighbour landscape, but evidently their understanding requires the use of theoretical and numerical tools, and presumably the development of new conceptual and practical approaches, of general relevance at the crossroads of mathematical logic, descriptive set theory, combinatorial game theory, game theory on graphs, etc.

 We complete this Section and the present work with a selected list of such open problems.

 \medskip

 \textbf{Completion of the dynamical landscape.} Construct explicit examples for the currently hypothetical trajectory types depicted in Figure \ref{fig:game_topology} and not detected (yet?), specifically: merging trajectories that proceed to a finite termination (collisions), and loops accessible via a single unique transient history (non-branching entry).

\medskip

\textbf{Taxonomy of initial decks.} Partition the set of balanced initial configurations $\mathcal{S}_{\textrm{start}}$ based on their asymptotic behaviour (terminating vs.\ infinite). Further classify decks into equivalence classes belonging to the same dynamical orbit, identifying the basins of attraction for specific loops.

\medskip

\textbf{Taxonomy of final winner's decks.} Analyse the statistical properties of terminal configurations $s_{\texttt{end}}$, focussing on meaningful indicators, be they the total number of cards, or the distribution, mean separation, and position entropy of special cards, or others. Determine if traces of the match's complexity or duration survive in the final card arrangement.

\medskip

\textbf{Identification of meaningful metrics in the state space.} Define distance functions on $\mathcal{S}$ capable of quantifying the convergence of trajectories or the proximity to termination. Address the current lack of predictive indicators for match length by identifying structural signatures that distinguish ultra-long transients from true periodic capture.

\medskip

\textbf{Understanding of the problem of backward determinism.} Rigorously characterise the failure of injectivity for the trick function $\mathfrak{F}$ on the reachable set, as discussed in Section \ref{sec:fwd-bkw-determinism}. Specifically, determine whether injectivity is lost only when a trajectory enters a loop, or also when two distinct trajectories merge and proceed to termination.

\medskip

\textbf{Completing the landscape of infinite loops and infinite matches.} Explore the boundaries of the `existence window' for loops and the conjecture that infinite games are confined to a specific intermediate regime of the special-to-total card ratio $4\mathfrak{R}/N$.

%

\begin{thebibliography}{10}

\bibitem{Berlekamp-Conway-Guy-2004}
{\sc E.~R. Berlekamp, J.~H. Conway, and R.~K. Guy}, {\em {Winning Ways for Your
  Mathematical Plays. Volume 4.}}, A K Peters, Wellesley MA, second~ed., 2004.

\bibitem{Casella2024}
{\sc B.~Casella, P.~M. Anderson, M.~Kleber, R.~P. Mann, R.~Nessler,
  W.~Rucklidge, S.~G. Williams, and N.~Wu}, {\em {A Non-Terminating Game of
  Beggar-My-Neighbor}}, The American Mathematical Monthly, 132 (2025),
  pp.~978--994.

\bibitem{Beggar-Collins-records}
{\sc T.~Collins}, {\em {Long playing card configurations in Beggar My
  Neighbor}}, http://www.tkcs-collins.com/truman/byn/byn.shtml.

\bibitem{Beggar-Gentilini}
{\sc A.~Gentilini}, {\em {Pelagaletto2}},
  https://github.com/alessandro-gentilini/pelagaletto2.

\bibitem{Lakshtanov-Aleksenko-Beggar2013}
{\sc E.~L. Lakshtanov and A.~I. Aleksenko}, {\em {Finiteness in the
  Beggar-My-Neighbor card game}}, Large Systems, 49 (2013), pp.~163--166.

\bibitem{Beggar-Mann-records}
{\sc R.~Mann}, {\em {Known historical Beggar-My-Neigbour records}},
  https://www.richardpmann.com/beggar-my-neighbour-records.html.

\bibitem{Beggar-Mayer-codes-records}
{\sc M.~Mayer}, {\em {Beggarmypython}},
  https://github.com/matthewmayer/beggarmypython.

\bibitem{Paulhus_Beggar1999}
{\sc M.~M. Paulhus}, {\em {Beggar My Neighbour}}, The American Mathematical
  Monthly, 106 (1999), pp.~162--65.

\bibitem{Spivey2010}
{\sc M.~Z. Spivey}, {\em {Cycles in war.}}, Integers, 10 (2010), pp.~747--764.

\bibitem{Beggar-Tristan-records}
{\sc F.-R. Tristan}, {\em {Beggar-My-Neigbour}},
  https://github.com/tristan-f-r/beggar-my-neighbour.

\bibitem{Beggar-Zanotto-codes}
{\sc R.~Zanotto}, {\em {Cavacamisa}},
  https://github.com/drago-96/cavacamisa/tree/master.

\end{thebibliography}

\def\cprime{$'$}

\end{document}